\documentclass[11pt]{amsart}
\usepackage{amssymb,amsthm,amsmath,epsfig,latexsym}
\usepackage{calc,times,verbatim}
\usepackage{graphicx}
\usepackage{pgf,tikz}
\usepackage{mathrsfs}
\usetikzlibrary{arrows}
\usepackage{geometry} 
\begin{document}

\newcommand{\mmbox}[1]{\mbox{${#1}$}}
\newcommand{\proj}[1]{\mmbox{{\mathbb P}^{#1}}}
\newcommand{\Cr}{C^r(\Delta)}
\newcommand{\CR}{C^r(\hat\Delta)}
\newcommand{\affine}[1]{\mmbox{{\mathbb A}^{#1}}}
\newcommand{\Ann}[1]{\mmbox{{\rm Ann}({#1})}}
\newcommand{\caps}[3]{\mmbox{{#1}_{#2} \cap \ldots \cap {#1}_{#3}}}
\newcommand{\Proj}{{\mathbb P}}
\newcommand{\N}{{\mathbb N}}
\newcommand{\Z}{{\mathbb Z}}
\newcommand{\R}{{\mathbb R}}
\newcommand{\A}{{\mathcal{A}}}
\newcommand{\Tor}{\mathop{\rm Tor}\nolimits}
\newcommand{\Ext}{\mathop{\rm Ext}\nolimits}
\newcommand{\Hom}{\mathop{\rm Hom}\nolimits}
\newcommand{\im}{\mathop{\rm Im}\nolimits}
\newcommand{\rank}{\mathop{\rm rank}\nolimits}
\newcommand{\supp}{\mathop{\rm supp}\nolimits}
\newcommand{\arrow}[1]{\stackrel{#1}{\longrightarrow}}
\newcommand{\CB}{Cayley-Bacharach}
\newcommand{\coker}{\mathop{\rm coker}\nolimits}
\sloppy
\theoremstyle{plain}

\newtheorem*{thm*}{Theorem}
\newtheorem{defn0}{Definition}[section]
\newtheorem{prop0}[defn0]{Proposition}
\newtheorem{quest0}[defn0]{Question}
\newtheorem{thm0}[defn0]{Theorem}
\newtheorem{lem0}[defn0]{Lemma}
\newtheorem{corollary0}[defn0]{Corollary}
\newtheorem{example0}[defn0]{Example}
\newtheorem{remark0}[defn0]{Remark}
\newtheorem{conj0}[defn0]{Conjecture}

\newenvironment{defn}{\begin{defn0}}{\end{defn0}}
\newenvironment{conj}{\begin{conj0}}{\end{conj0}}
\newenvironment{prop}{\begin{prop0}}{\end{prop0}}
\newenvironment{quest}{\begin{quest0}}{\end{quest0}}
\newenvironment{thm}{\begin{thm0}}{\end{thm0}}
\newenvironment{lem}{\begin{lem0}}{\end{lem0}}
\newenvironment{cor}{\begin{corollary0}}{\end{corollary0}}
\newenvironment{exm}{\begin{example0}\rm}{\end{example0}}
\newenvironment{rem}{\begin{remark0}\rm}{\end{remark0}}

\newcommand{\defref}[1]{Definition~\ref{#1}}
\newcommand{\conjref}[1]{Conjecture~\ref{#1}}
\newcommand{\propref}[1]{Proposition~\ref{#1}}
\newcommand{\thmref}[1]{Theorem~\ref{#1}}
\newcommand{\lemref}[1]{Lemma~\ref{#1}}
\newcommand{\corref}[1]{Corollary~\ref{#1}}
\newcommand{\exref}[1]{Example~\ref{#1}}
\newcommand{\secref}[1]{Section~\ref{#1}}
\newcommand{\remref}[1]{Remark~\ref{#1}}
\newcommand{\questref}[1]{Question~\ref{#1}}

\newcommand{\std}{Gr\"{o}bner}
\newcommand{\jq}{J_{Q}}

\title{Betti numbers of fold products of linear forms}
\author{Ricardo Burity and \c{S}tefan Toh\v{a}neanu}

 \keywords{Betti numbers, free resolution, fold products. \\
		\indent The first author was partially supported by a grant from CNPq (Brazil)(408698/2023-3).\\
\indent Burity's address: Departamento de Matemática, Universidade Federal da Paraiba, 58051-900, J. Pessoa, PB, , Brazil,
		Email: ricardo@mat.ufpb.br.\\
\indent Tohaneanu's address: Department of Mathematics, University of Idaho, Moscow, Idaho 83844-1103, USA, Email: tohaneanu@uidaho.edu.}

\begin{abstract}
The ideals generated by fold products of linear forms are generalizations of powers of defining ideals of star configurations, or of Veronese type ideals, and in this paper we study their Betti numbers. In earlier work, the authors together with Yu Xie showed that these ideals have linear graded free resolution, and in this paper we take the study of the homological information about these ideals to the next natural level.
\end{abstract}

\maketitle

\section{Introduction}

Let $R:=\mathbb K[x_1,\ldots,x_k]$ be the ring of polynomials with coefficients in a field $\mathbb K$, with the standard grading. Denote ${\frak m}:=\langle x_1,\ldots,x_k\rangle$ to be the irrelevant maximal ideal of $R$. Let $\ell_1,\ldots,\ell_n$ be linear forms in $R$, some possibly proportional, and denote this collection by $\Sigma=(\ell_1,\ldots,\ell_n)\subset R$; we need a notation to take into account the fact that some of these linear forms are proportional. For $\ell\in \Sigma$, by $\Sigma\setminus\{\ell\}$ we will understand the collection of linear forms of $\Sigma$ from which $\ell$ has been removed. Also, we denote $|\Sigma|=n$, and ${\rm rk}(\Sigma):={\rm ht}(\langle \ell_1,\ldots,\ell_n\rangle)$.

Let $1\leq a\leq n$ be an integer and define {\em the ideal generated by $a$-fold products of $\Sigma$} to be the ideal of $R$ $$I_a(\Sigma):=\langle \{\ell_{i_1}\cdots\ell_{i_a}|1\leq i_1<\cdots<i_a\leq n\}\rangle.$$ We also make the convention $I_0(\Sigma):=R$, and $I_b(\Sigma)=0$, for all $b>n$. Also, if $\Sigma=\emptyset$, $I_a(\Sigma)=0$, for any $a\geq 1$.

When the linear forms $\ell_1,\ldots,\ell_n$ define a hyperplane arrangement $\A\subset\mathbb P^{k-1}$ (so for all $i\neq j$, $\gcd(\ell_i,\ell_j)=1$), we will denote $\Sigma=\A=\{\ell_1,\ldots,\ell_n\}$.

\medskip

A homogeneous ideal $I\subset R$ generated in degree $a$ it is said to have {\em linear (minimal) graded free resolution}, if one has the graded free resolution $$0\rightarrow R^{b_{t}}(-(a+t-1))\rightarrow\cdots\rightarrow R^{b_2}(-(a+1))\rightarrow R^{b_1}(-a)\rightarrow I\rightarrow 0,$$ for some positive integer $t$. The integers $b_j\geq 1$ are called the {\em Betti numbers} of $I$ (or of $R/I$).\footnote{In general, for any homogeneous ideal $I$, the Betti numbers of $I$ are the ranks of the free modules in a minimal free resolution of $I$. But because the minimal free resolution is most of the time not pure, the attention is directed towards the ranks of each graded piece in the free modules, ranks that are called {\em graded Betti numbers}.} By convention, the zero ideal has linear graded free resolution. Also we say that $R/I$ has linear graded free resolution if and only if $I$ has linear graded free resolution.

\begin{rem}\label{basic} For any ideal $I\subset R=\mathbb K[x_1,\ldots,x_k]$, the following are true for the Betti numbers of $I$.

\begin{itemize}
  \item[(i)] If $b_s=0$ for some $s\geq 1$, then $b_t=0$ for all $t\geq s$.
  \item[(ii)] If $p:={\rm pdim}_R(R/I)$ is the projective dimension, then $b_{p+1}=0$ and $b_p\neq 0$.
  \item[(iii)] Hilbert Syzygy Theorem: $k\geq {\rm pdim}_R(R/I)$, and so $b_{k+1}=0$.
\end{itemize}

We will use these basic results about vanishing of some Betti numbers without explicitly calling upon this remark.
\end{rem}

\cite[Conjecture 1]{AnGaTo} states that for any collection of linear forms $\Sigma$, and any $1\leq a\leq|\Sigma|$, the ideals $I_a(\Sigma)$ have linear graded free resolution. In \cite[Theorem 2.2]{BuToXi} this conjecture was proved, and so we have:

$$0\longrightarrow R^{b_{t}(a,\Sigma)}(-(a+t-1))\longrightarrow \cdots \longrightarrow R^{b_2(a,\Sigma)}(-(a+1))\longrightarrow R^{b_1(a,\Sigma)}(-a)\longrightarrow I_a(\Sigma)\longrightarrow 0,$$ for some integer $t\geq 1$.

The goal is to find formulas for all $b_i(a,\Sigma)$.

Finding the complete description of the Betti numbers of ideals with linear free resolution is an interesting by itself; see, for example, \cite{EiGo, HeKu, Ni}, and the citations therein. In the Cohen-Macauly case, it is ``rather easy'' to obtain this information; see \cite[Proposition 1.7]{EiGo}. The coordinate ring of a star configuration (see, \cite{GeHaMi}) is the quotient of $R$ by an ideal generated by fold products of linear forms, and it is also Cohen-Macaulay. Though \cite[Proposition 2.1]{GeHaMi} doesn't appeal to Eisenbud-Goto result to obtain all the homological information, it bares the question what happens when we quotient $R$ by the defining ideal of a ``generalized star configuration'',\footnote{Some authors understand by ``generalized star configuration'' as the analog of a star configuration but where the linear forms defining the skeleton are replaced by sufficiently generic forms of the same degree $d\geq 2$.}  meaning any ideal generated by fold products of linear forms. Sections 2.2 and 3.2 of \cite{Ni} provide such computations for ideals generated by fold products of linear forms. Also \cite[Theorem 2.4]{To3} and \cite[Lemma 3.2]{Sc} present the Betti numbers for the ideal generated by the $(n-2)$, and respectively the $(n-1)$, fold products of $n$ nonproportional linear forms. Also, \cite[Proposition 2.10]{AnGaTo} presents the combinatorial formula of the first Betti number for any ideal generated by fold products of linear forms.

In Section \ref{sec_known} we present all the known results, some of which we mentioned above, and immediate consequences of some known results and methods. We summarize these in Theorem \ref{thm_summary} at the end of the section. Note that the results are strongly dependent on the generalized Hamming weights of the linear code with generator matrix the coefficients matrix of the linear forms. In Section \ref{sec_DelContr} we present a recursive way to compute the Betti numbers via the common technique of deletion-contraction; see Theorem \ref{thm_recursion}. In Section \ref{sec_k3}, we apply this theorem for some special case of collections of linear forms in $\mathbb K[x_1,x_2,x_3]$: see Corollaries \ref{cor_k3generic} and \ref{cor_dplus1}. In the following Section \ref{sec_comb} we present another method to obtain the Betti numbers via some combinatorics and linear algebra (Lemma \ref{lem_HF} via Theorem \ref{thm_kernel}). As reflected in Example \ref{exm_kernelPsi} and Remark \ref{rem_revisit} this method is rather difficult to apply in general, and so it is our belief that Theorem \ref{thm_recursion} is still the most computationally efficient way of finding the Betti numbers (see also the discussions in Remark \ref{rem_Veronese} and after it). Not only that, but since the recursion boils down to finding Betti numbers when $k=2$, and since deletion-contraction depends only on the combinatorics, we can conclude that all Betti numbers of any ideal generated by fold products of linear forms are combinatorially determined.

\section{Known Betti numbers formulas and calculations}\label{sec_known}

Let $\Sigma=(\ell_1,\ldots,\ell_n)$ be a collection of linear forms in $R:=\mathbb K[x_1,\ldots,x_k]$. Let $1\leq a\leq n$ be an integer. We first present a section of known results about the Betti numbers concerning various cases of $I_a(\Sigma)$. After a change of variables, we can assume that $n\geq k$ and that $\langle \ell_1,\ldots,\ell_n\rangle=\langle x_1,\ldots,x_k\rangle=:\frak m$.

\medskip

\subsection{Minimum distance.}\label{mindist} One can think of the linear forms in $\Sigma$ as the dual forms of some points in $\mathbb P^{k-1}$, counted with repetitions. Fixing representatives for their homogeneous coordinates, these points can be placed as the columns of an $k\times n$ matrix, $G_{\Sigma}$, in any order we wish. After a change of variables, we can assume that $G_{\Sigma}$ has full rank equal to $k$.

The matrix $G_{\Sigma}$ can be thought as a generator matrix of a linear code $\mathcal C_{\Sigma}$. Changing the representatives or the order of the columns will produce a linear code ``equivalent'' to $\mathcal C_{\Sigma}$; the parameters of this ``new'' code are exactly the same as of the old one.

{\em The minimum distance} of $\mathcal C_{\Sigma}$ is the minimum number of nonzero entries in a nonzero element of the row-span of $G_{\Sigma}$. We denote this invariant with $d(\Sigma)$. By \cite[Theorem 3.1]{To}, if $1\leq a\leq d(\Sigma)$, then $I_a(\Sigma)=\frak m^a$. Then, by \cite[Proposition 1.7(c)]{EiGo}, the nonzero Betti numbers of $I_a(\Sigma)$ are $$b_i(a,\Sigma)={{k+a-1}\choose{a+i-1}}\cdot{{a+i-2}\choose{a-1}}, i=1,\ldots,k.$$

Let $1\leq r\leq k$ be an integer. The {\em $r$-th generalized Hamming weight} is an invariant of $\mathcal C_{\Sigma}$ that can be calculated in the following way: $d_r(\Sigma)$ is $n$ minus the maximum number of columns of $G_{\Sigma}$ that span a $(k-r)$-dimensional vector space (see \cite[Proposition 2.8]{To0}). If $r=1$, one recovers the minimum distance $d(\Sigma)$ from above; if $r=k$, then $d_k(\Sigma)=n$, since $G_{\Sigma}$ doesn't have any zero columns.

By \cite[Proposition 2.2]{AnGaTo}, the generalized Hamming weights help us figure out the height of $I_a(\Sigma)$. With the convention that $d_0(\Sigma)=1$, for $r=0,\ldots,k-1$, if $d_r(\Sigma)<a\leq d_{r+1}(\Sigma)$, then ${\rm ht}(I_a(\Sigma))=k-r$.

\begin{rem}\label{rem_height1} An interesting case would be to understand when ${\rm ht}(I_a(\Sigma))=1$. This happens whenever $d_{k-1}(\Sigma)+1\leq a\leq n$, and we know that $d_{k-1}(\Sigma)=n-w$, where $w$ is the maximum number of proportional linear forms in $\Sigma$.

Suppose $\Sigma=(\underbrace{\ell_1,\ldots,\ell_1}_{m_1},\ldots,\underbrace{\ell_t,\ldots,\ell_t}_{m_t})$, with $t\geq 1$, $\gcd(\ell_i,\ell_j)=1$, $w=m_1\geq\cdots\geq m_t\geq 1$, $m_1+\cdots+m_t=n$, and $\langle \ell_1,\ldots,\ell_t\rangle=\frak m$.

If $w=1$, then ${\rm ht}(I_a(\Sigma))=1$ iff $a=n$. In this case, the only nonzero Betti number is $$b_1(n,\Sigma)=1.$$

Suppose $w\geq 2$. Suppose that for $i=1,\ldots,s$, $e_i:=m_i+a-n\geq 1$, and $n-a\geq m_{s+1}\geq\cdots\geq m_t$. It is not hard to see that $$I_a(\Sigma)=\ell_1^{e_1}\cdots\ell_s^{e_s} I_{a-(e_1+\cdots +e_s)}(\Sigma'),$$ where $\Sigma'=(\underbrace{\ell_1,\ldots,\ell_1}_{n-a},\ldots,\underbrace{\ell_s,\ldots,\ell_s}_{n-a}, \underbrace{\ell_{s+1},\ldots,\ell_{s+1}}_{m_{s+1}},\underbrace{\ell_t,\ldots,\ell_t}_{m_t})$, since $m_i-e_i=n-a$, for $i=1,\ldots,s$.

But, since $a-(e_1+\cdots+e_s)=s(n-a)+m_{s+1}+\cdots+m_t-(n-a)$, and since $\max\{n-a,m_{s+1},\ldots,m_t\}=n-a$, we have that $a-(e_1+\cdots+e_s)=d_{k-1}(\Sigma')$, and therefore, ${\rm ht}(I_{a-(e_1+\cdots +e_s)}(\Sigma'))=2$.

Nonetheless, in terms of Betti numbers we have $$b_i(a,\Sigma)=b_i(a-(e_1+\cdots+e_s),\Sigma'), i=1,\ldots,k.$$
\end{rem}

\begin{lem}\label{rank_two} Let $\Sigma=(\underbrace{\ell_1,\ldots,\ell_1}_{m_1},\ldots,\underbrace{\ell_t,\ldots,\ell_t}_{m_t})\subset R:=\mathbb K[x_1,x_2]$, with $t\geq 1$, and $\gcd(\ell_i,\ell_j)=1$. Let $1\leq a\leq n:=m_1+\cdots+m_t$, and for $i=1,\ldots, t$, let $e_i:=\max\{m_i+a-n,0\}$. Let $e:=\max\{a-(e_1+\cdots+e_t),0\}$. Then,
$$I_a(\Sigma)=\langle\ell_1\rangle^{e_1}\cap\cdots\cap\langle\ell_t\rangle^{e_t}\cap \langle x_1,x_2\rangle^a=\ell_1^{e_1}\cdots\ell_t^{e_t}\cdot\langle x_1,x_2\rangle^e,$$ and so the nonzero Betti numbers are $$b_1(a,\Sigma)=e+1 \mbox{ and } b_2(a,\Sigma)=e.$$
\end{lem}
\begin{proof} The proofs follows immediately from the argument in Remark \ref{rem_height1}, which was inspired by \cite[Proposition 2.2]{BuToXi} and \cite[Theorem 2.2]{To3} for the case $k=2$.

Note that in this formula we included also the case when $a\leq n-\max\{m_1,\ldots,m_t\}$, which is $d(\Sigma)$. In this case, $e_1=\cdots=e_t=0$ and $e=a$; exactly the statement of \cite[Theorem 3.1]{To}.
\end{proof}
\medskip

\subsection{The Cohen-Macaulay case.}\label{CM} When $1\leq a\leq d(\Sigma)$, observe that $R/I_a(\Sigma)$ is (arithmetically) Cohen-Macaulay, and therefore we were able to apply the Eisenbud-Goto result to obtain the corresponding formulas for the Betti numbers. The same result will be applied in other instances when $R/I_a(\Sigma)$ is Cohen-Macaulay; i.e, ${\rm pdim}_R(R/I_a(\Sigma))=p={\rm ht}(I_a(\Sigma))$.

By \cite[Theorem 3.38]{To0},$${\rm pdim}_R(R/I_a(\Sigma))=\min\{k,n-a+1\}.$$ Therefore, if ${\rm ht}(I_a(\Sigma))<k$, then $R/I_a(\Sigma)$ is Cohen-Macaulay iff ${\rm ht}(I_a(\Sigma))=n-a+1$; this implies that $n-a+1<k$ and so $n-k+1<a\leq n$. (see also \cite[Corollaries 3.13 and 3.39]{To0}).

Let $h\leq k$ be a positive integer. $\Sigma$ is called {\em $(h-1)$-generic} if any $h$ linear forms of $\Sigma$ are linearly independent. ``$(k-1)$-generic'' is the definition of what is commonly known as ``generic hyperplane arrangement''. Also, ``1-generic'' means that no two elements of $\Sigma$ are proportional, and ``0-generic'' means that $0\notin\Sigma$. Of course, if $\Sigma$ is $(h-1)$-generic, then it is $(h-u)$-generic for any $u\in\{1,\ldots,h\}$.

\begin{lem}\label{lem_generic} If $\Sigma$ is $(n-a)$-generic, where $n-k+1\leq a\leq n$, then the nonzero Betti numbers $I_a(\Sigma)$ are
$$b_i(a,\Sigma)={{n}\choose{i+a-1}}\cdot{{i+a-2}\choose{a-1}}, i=1,\ldots,n-a+1.$$
\end{lem}
\begin{proof} If $n-k+1\leq a\leq n$, then $1\leq h:=n-a+1\leq k$, and by \cite[Lemma 3.41]{To0}, $\Sigma$ is $(n-a)$-generic iff ${\rm ht}(I_a(\Sigma))=n-a+1$. From above, this is equivalent to $R/I_a(\Sigma)$ being Cohen-Macaulay. Therefore, the claim follows immediately from \cite[Proposition 1.7(c)]{EiGo}.
\end{proof}

When $a=n-c+1$, for some integer $1\leq c\leq k-1$, and $\Sigma$ is $(k-1)$-generic (which implies that it is also $(n-a)$-generic, as $n-a=c-1\leq k-2$), then $I_a(\Sigma)$ defines a ``codimension $c$ star configuration'' (see \cite[Lemma 3.1]{ToXi}), with Betti numbers $\displaystyle b_i={{n}\choose{c-i}}\cdot{{n-c+i-1}\choose{i-1}}, i=1,\ldots,c$. Note that, after some computations with binomial coefficients, this matches with what one will obtain when using the $h$-vector presented in \cite[Proposition 2.9 (2)]{GeHaMi}.

\medskip

\subsection{The Tutte polynomial.}\label{tutte} By definition, {\em the Tutte polynomial} of a matroid ${\rm M}$ of rank $k$ on the ground set $[n]=\{1,\ldots,n\}$ with rank function $r$ is
$$
T_{{\rm M}}(x,y)=\sum_{I\subseteq [n]}(x-1)^{k-r(I)}(y-1)^{|I|-r(I)}.
$$ In our case, the matroid ${\rm M}$ is the matroid of the matrix $G_{\Sigma}$ from above (the matrix of coefficients of the linear forms of $\Sigma$). If $I\subseteq [n]$, then our $r(I)$ is the rank of the submatrix of $G_{\Sigma}$ with columns indexed by $I$.

By \cite[Proposition 2.10]{AnGaTo}, if $1\leq a\leq n$, then the first Betti number equals
$$b_1(a,\Sigma)=\sum_{u=0}^{\min\{k,n-a\}}c_{k-u,n-a-u},$$ where $c_{i,j}$ is the coefficient of $x^iy^j$ in the expansion of $T_{{\rm M}}(x+1,y)$.

\begin{exm}\label{example_tutte} Consider $\Sigma=(x_1, x_1, x_2, x_3, x_1-x_3, x_2+x_3, x_1+2x_2+5x_3)\subset R:=\mathbb K[x_1,x_2,x_3]$. We have $k=3$, $n=6$ and $\displaystyle G_{\Sigma}=\left[\begin{array}{rrrrrrr}1&1&0&0&1&0&1\\
0&0&1&0&0&1&2\\0&0&0&1&-1&1&5\end{array}\right]$.

The (affine) picture is presented below where we marked by bold points the multiple intersection points. For simplicity, we label the lines by their defining linear forms.

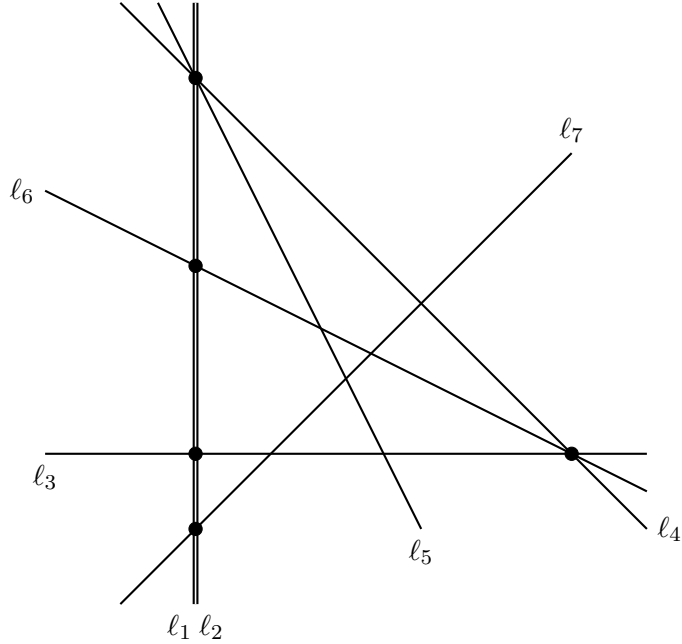
\begin{center}
\begin{tikzpicture}
\draw[thick] (-2,0) node[below]{$\ell_3$} -- (6,0);
\draw[thick, double] (0,-2) node[below]{$\ell_1$ $\ell_2$} -- (0,6);
\draw[thick] (-1,6) -- (6,-1) node[right]{$\ell_4$};
\draw[thick] (-2,3.5) node[left]{$\ell_6$} -- (6,-0.5);
\draw[thick] (-0.5,6) -- (3,-1) node[below]{$\ell_5$};
\draw[thick] (-1,-2) -- (5,4) node[above]{$\ell_7$};
\filldraw[black] (0,0) circle (2.5pt);
\filldraw[black] (0,2.5) circle (2.5pt);
\filldraw[black] (0,5) circle (2.5pt);
\filldraw[black] (5,0) circle (2.5pt);
\filldraw[black] (0,-1) circle (2.5pt);
\end{tikzpicture}
\end{center}

We can use \cite{GrSt} to obtain the Tutte polynomial of the matroid of this matrix, but we decided to give more details below.

\begin{itemize}
  \item[(a)] $I=\emptyset$ has $r(I)=0$ and $|I|=0$, so it contributes with the summand $(x-1)^3$.
  \item[(b)] $I=\{j\}, j=1,\ldots,7$ have $r(I)=1$ and $|I|=1$, so they contribute with the summand $7(x-1)^2$.
  \item[(c)] $I=\{1,2\}$ has $r(I)=1$ and $|I|=2$, so it contributes with the summand $(x-1)^2(y-1)$.
  \item[(d)] $I=\{j_1,j_2\}, j_1<j_2$ except for the one at (c) have $r(I)=2$ and $|I|=2$, so they contribute with the summand $({{7}\choose{2}}-1)(x-1)=20(x-1)$.
  \item[(e)] $I=\{1,2,j\}, j=3,\ldots,7$ and $\{1,4,5\}, \{2,4,5\}, \{3,4,6\}$ have $r(I)=2$ and $|I|=3$, so they contribute with the summand $8(x-1)(y-1)$.
  \item[(f)] $I=\{1,2,4,5\}$ has $r(I)=2$ and $|I|=4$, so it contributes with the summand $(x-1)(y-1)^2$.
  \item[(g)] $I=\{j_1,j_2,j_3\}, j_1<j_2<j_3$ except for the ones at (e) have $r(I)=3$ and $|I|=3$, so they contribute with the summand $({{7}\choose{3}}-8)=27$.
  \item[(h)] $I=\{j_1,j_2,j_3,j_4\},j_1<j_2<j_3<j_4$ except for the one at (f) have $r(I)=3$ and $|I|=4$, so they contribute with the summand $({{7}\choose{4}}-1)(y-1)=34(y-1)$.
  \item[(i)] $I=\{j_1,j_2,j_3,j_4,j_5\}, j_1<\cdots<j_5$ have $r(I)=3$ and $|I|=5$, so they contribute with the summand $21(y-1)^2$.
  \item[(j)] $I=\{j_1,\ldots,j_6\},j_1<\cdots<j_6$ have $r(I)=3$ and $|I|=6$, so they contribute with the summand $7(y-1)^3$.
  \item[(k)] $I=[7]$ has $r(I)=3$ and $|I|=7$, so it contributes with the summand $(y-1)^4$.
\end{itemize}

The Tutte polynomial evaluated at $(x+1,y)$ is

$$T_{{\rm M}}(x+1,y)=y^4+x^3+x^2y+xy^2+3y^3+6x^2+6xy+6y^2+13x+9y+8.$$

In this case, $d(\Sigma)=7-4=3$, so the first three values in the second column in the table below are immediate since $I_a(\Sigma)=\langle x_1,x_2,x_3\rangle^a, a=1,2,3$. Also, the last value in this column is also known since $I_7(\Sigma)=\langle x_1^2x_2x_3(x_1-x_3)(x_2+x_3)(x_1+2x_2+5x_3)\rangle$. The other values in this column are computed using \cite{GrSt}, and they match the values in the third column.

\medskip

\begin{center}
\begin{tabular}{|c|c|l|}
\hline
$a$ & $b_1(a,\Sigma)$ & formula from Tutte polynomial\\
\hline
$1$ & $3$ & $c_{3,6}+c_{2,5}+c_{1,4}+c_{0,3}=0+0+0+3=3$\\
\hline
$2$ & $6$ & $c_{3,5}+c_{2,4}+c_{1,3}+c_{0,2}=0+0+0+6=6$\\
\hline
$3$ & $10$ & $c_{3,4}+c_{2,3}+c_{1,2}+c_{0,1}=0+0+1+9=10$\\
\hline
$4$ & $14$ & $c_{3,3}+c_{2,2}+c_{1,1}+c_{0,0}=0+0+6+8=14$\\
\hline
$5$ & $14$ & $c_{3,2}+c_{2,1}+c_{1,0}=0+1+13=14$\\
\hline
$6$ & $6$ & $c_{3,1}+c_{2,0}=0+6=6$\\
\hline
$7$ & $1$ & $c_{3,0}=1$\\
\hline
\end{tabular}
\end{center}
\end{exm}

\medskip

\subsection{Herzog-K\"{u}hl Equations.}\label{HerzogKuhl} An ideal $I$ is said to have {\em linear powers} if $I^u$ has linear graded free resolution for all $u\geq 1$.

For any $1\leq a\leq n$, our ideals $I_a(\Sigma)$ have linear powers: if $u\geq 1$, then $(I_a(\Sigma))^u=I_{ua}(\Sigma(u))$, where $\Sigma(u)=(\underbrace{\ell_1,\ldots,\ell_1}_{u},\ldots,\underbrace{\ell_n,\ldots,\ell_n}_{u})$, and the later has linear graded free resolution by \cite[Theorem 2.2]{BuToXi} (see also the proof of Lemma 3.1 in \cite{BuToXi}).

With this in mind, if $\delta:=k-{\rm ht}(I_a(\Sigma))$, we have the following {\em $k-\delta$ Herzog-K\"{u}hl equations} (see \cite{HeKu}):

\begin{eqnarray}\label{HHequations}
\sum_{i=1}^k(-1)^ib_i(a,\Sigma)=-1\nonumber\\
\sum_{i=1}^k(-1)^i(a+i-1)b_i(a,\Sigma)=0\nonumber\\
\sum_{i=1}^k(-1)^i(a+i-1)(a+i-2)b_i(a,\Sigma)=0\nonumber\\
\vdots\nonumber\\
\sum_{i=1}^k(-1)^i(a+i-1)\cdots(a+i-(k-\delta-1))b_i(a,\Sigma)=0.\nonumber
\end{eqnarray}

\begin{rem}\label{rem_rank3} If $k=3$, then ${\rm ht}(I_a(\Sigma))$ can be $1, 2$ or $3$.

\begin{itemize}
  \item[(i)] When ${\rm ht}(I_a(\Sigma))=3$, then, by \cite{To}, $I_a(\Sigma)=\frak m^a$, a case we discussed in Section \ref{mindist} above.
  \item[(ii)] If ${\rm ht}(I_a(\Sigma))=2$, then we have two Herzog-K\"{u}hl equations in three Betti numbers: $$b_1(a,\Sigma)-b_2(a,\Sigma)+b_3(a,\Sigma)=1$$ and $$ab_1(a,\Sigma)-(a+1)b_2(a,\Sigma)+(a+2)b_3(a,\Sigma)=0.$$ Since $b_1(a,\Sigma)$ can be computed from Section \ref{tutte}, we have the complete description of the Betti numbers:
      $$b_2(a,\Sigma)=2b_1(a,\Sigma)-a-2 \mbox{ and } b_3(a,\Sigma)=b_1(a,\Sigma)-a-1.$$
  \item[(iii)] Suppose ${\rm ht}(I_a(\Sigma))=1$. Then, from Section \ref{mindist}, $$b_i(a,\Sigma)=b_i(a-(e_1+\cdots+e_s),\Sigma'), i=1,2,3,$$ where ${\rm ht}(I_{a-(e_1+\cdots+e_s)}(\Sigma'))=2$, hence its Betti numbers can be determined from (ii) above.
\end{itemize}
\end{rem}

\begin{exm}\label{exm_betti_tutte} In Example \ref{tutte}, the cases not immediate are when $a=4,5,$ and $6$. Using \cite{GrSt} for height calculations, we have:

\begin{itemize}
  \item $\boxed{a=4}$ In this case, ${\rm ht}(I_4(\Sigma))=2$, and since $b_1(4,\Sigma)=14$, we have
  $$b_2(4,\Sigma)=2\cdot 14-4-2=22 \mbox{ and } b_3(4,\Sigma)=14-4-1=9.$$
  \item $\boxed{a=5}$ In this case, ${\rm ht}(I_5(\Sigma))=2$, and since $b_1(5,\Sigma)=14$, we have
  $$b_2(5,\Sigma)=2\cdot 14-5-2=21 \mbox{ and } b_3(5,\Sigma)=14-5-1=8.$$
  \item $\boxed{a=6}$ In this case, ${\rm ht}(I_6(\Sigma))=1$. Then, by Remark \ref{rem_height1} with $e_1=1$ and $e_i=0, i=2,\ldots,5$, giving that $b_i(6,\Sigma)=b_i(5,(x_1,x_2,x_3,x_1-x_3,x_2+x_3,x_1+2x_2+5x_3)),i=1,2,3$. From Section \ref{hyperplanes} below, part (ii), we then have
      $$b_1(6,\Sigma)=6, b_2(6,\Sigma)=5, b_3(6,\Sigma)=0.$$
\end{itemize}
\end{exm}

\begin{rem}\label{rem_points} More generally than part (ii) in Remark \ref{rem_rank3}, is another situation when the Herzog-K\"{u}hl equations can determine completely the Betti numbers (in terms of the first Betti number): $d_1(\Sigma)+1\leq a\leq d_2(\Sigma)$. This is because if $a$ is in this range, then ${\rm ht}(I_a(\Sigma))=k-1$. We then have $k-1$ equations in $k$ parameters (i.e., the Betti numbers), and \cite[Theorem 3]{Ni} gives the linear algebra calculation:

$$b_i(a,\Sigma)=\sum_{j=0}^{k-2}{{j}\choose{i-2}}\left(b_1(a,\Sigma)-{{a+j}\choose{j}}\right),\, i=2,\ldots,k.$$
\end{rem}

\medskip

\subsection{Hyperplane arrangements.}\label{hyperplanes} Suppose $\Sigma=(\ell_1,\ldots,\ell_n)\subset R:=\mathbb K[x_1,\ldots,x_k]$ defines an essential (central) hyperplane arrangement in $\mathbb K^k$. This means that $\ell_i$'s are linear forms, with $\gcd(\ell_i,\ell_j)=1, i\neq j$, and $\langle \ell_1,\ldots,\ell_n\rangle=\langle x_1,\ldots,x_k\rangle=:\frak m$.\footnote{As we said, we could have denoted this by $\Sigma=\{\ell_1,\ldots,\ell_n\}$, with ${\rm rk}(\Sigma)=k$.}

\medskip

(i) We already saw that if $a=n$, then $$b_1(n,\Sigma)=1 \mbox{ and } b_2(n,\Sigma)=0,$$ as $I_n(\Sigma)=\langle \ell_1\cdots\ell_n\rangle$, a principal ideal.

\medskip

(ii) If $a=n-1$, then $I_{n-1}(\Sigma)=\langle \ell_2\cdots\ell_n,\ldots,\ell_1\cdots\ell_{n-1}\rangle$. In the case of hyperplane arrangements, ${\rm ht}(I_{n-1}(\Sigma))=2$, and the Betti numbers are $$b_1(n-1,\Sigma)=n, b_2(n-1,\Sigma)=n-1, \mbox{ and } b_3(n-1,\Sigma)=0;$$ see, for example, \cite[Lemma 3.2]{Sc}.

\medskip

(iii) If $a=n-2$, then, by \cite[Theorem 2.4]{To3}, we have
$$b_1(n-2,\Sigma)=\alpha-\beta, b_2(n-2,\Sigma)=2\alpha-n-2\beta, b_3(n-2,\Sigma)=\alpha-n-\beta+1, \mbox{ and } b_4(n-2,\Sigma)=0,$$ where $\displaystyle \alpha:={{n}\choose{2}}$ and $\displaystyle \beta:=\sum_{X\in L_2(\Sigma)}{{|\Sigma_X|-1}\choose{2}}$. Here, $X\in L_2(\Sigma)$ means that $X=V(\ell_i,\ell_j)$ for some $i\neq j$, and $\Sigma_X=\{\ell\in\Sigma|\ell\in\langle \ell_i,\ell_j\rangle\}$.

One should observe that we obtained the same values for the first three Betti numbers as in part (ii) of Remark \ref{rem_rank3}.

\medskip

\subsection{Summary of known results.} We end this section with a result where we summarize the known results we discussed above.

\begin{thm}\label{thm_summary} Let $\Sigma=(\underbrace{\ell_1,\ldots,\ell_1}_{m_1},\ldots,\underbrace{\ell_t,\ldots,\ell_t}_{m_t})\subset R:=\mathbb K[x_1,\ldots,x_k]$ be a collection of linear forms, with $t\geq 1$, $\gcd(\ell_i,\ell_j)=1, i\neq j$, $m_1\geq\cdots\geq m_t\geq 1$, $m_1+\cdots+m_t=n$, and $\langle \ell_1,\ldots,\ell_t\rangle=\frak m$. For $1\leq r\leq k$, an integer, let $d_r(\Sigma)$ be the $r$-th generalized Hamming weight of $\Sigma$. Let $1\leq a\leq n$ be an integer, and for $i=1,\ldots,k$ let $b_i(a,\Sigma)$ denote the $i$-th Betti number of $I_a(\Sigma)$. Also, for $i=1,\ldots,t$, let $e_i:=\max\{m_i+a-n,0\}$ and $e:=\max\{a-(e_1+\cdots+e_t),0\}$. Then the following happen:
\begin{itemize}
  \item[(0)] If $G_{\Sigma}$ is the matrix of coefficients of $\Sigma$, then for $1\leq r\leq k$, $d_{r}(\Sigma)=n-M_{k-r}$, where $M_{k-r}$ is the maximum number of columns of $G_{\Sigma}$ that span a $(k-r)$-dimensional vector space.
  \item[(1)] If $T_{\rm M}(x,y)$ is the Tutte polynomials of the matroid ${\rm M}$ of $G_{\Sigma}$, then $$b_1(a,\Sigma)=\sum_{u=0}^{\min\{k,n-a\}}c_{k-u,n-a-u},$$ where $c_{i,j}$ is the coefficient of $x^iy^j$ in the expansion of $T_{{\rm M}}(x+1,y)$.
  \item[(2)] If $\boxed{1\leq a\leq d_1(\Sigma)}$, then $$b_i(a,\Sigma)={{k+a-1}\choose{a+i-1}}\cdot{{a+i-2}\choose{a-1}}, i=1,\ldots,k.$$
  \item[(3)] If $\boxed{d_1(\Sigma)+1\leq a\leq d_2(\Sigma)}$, then $$b_i(a,\Sigma)=\sum_{j=0}^{k-2}{{j}\choose{i-2}}\left(b_1(a,\Sigma)-{{a+j}\choose{j}}\right),\, i=2,\ldots,k.$$
  \item[(4)] We have $d_{k-1}(\Sigma)=n-m_1$ and $d_k(\Sigma)=n$. \\If $\boxed{d_{k-1}(\Sigma)+1\leq a\leq d_k(\Sigma)}$, then ${\rm ht}(I_a(\Sigma))=1$.
  \begin{itemize}
    \item[(4.i)] If $\boxed{a=n}$, then $$b_1(n,\Sigma)=1 \mbox{ and } b_2(n,\Sigma)=0.$$
    \item[(4.ii)] If $m_1=1$, then ${\rm ht}(I_a(\Sigma))=1$ iff $a=n$. In this case, $$b_1(n,\Sigma)=1 \mbox{ and } b_2(n,\Sigma)=0.$$
    \item[(4.iii)] If $m_1\geq 2$, then $$b_i(a,\Sigma)=b_i(e,\Sigma'), i=1,\ldots,k,$$ where $\Sigma'=(\underbrace{\ell_1,\ldots,\ell_1}_{m_1-e_1},\ldots,\underbrace{\ell_t,\ldots,\ell_t}_{m_t-e_t})$; note that $d_{k-1}(\Sigma')=e$.
  \end{itemize}
  \item[(5)] If $\boxed{a=n-1}$, then $$b_1(n-1,\Sigma)=t, b_2(n-1,\Sigma)=t-1 \mbox{ and } b_3(n-1,\Sigma)=0.$$
  \item[(6)] If $\boxed{k=2}$, then
    \begin{itemize}
    \item[(6.i)] if $1\leq a\leq d_1(\Sigma)=n-m_1$, then $$b_1(a,\Sigma)=a+1, b_2(a,\Sigma)=a, \mbox{ and } b_3(a,\Sigma)=0.$$
    \item[(6.ii)] if $d_1(\Sigma)+1\leq a\leq d_2(\Sigma)=n$, then $$b_1(a,\Sigma)=e+1, b_2(a,\Sigma)=e, \mbox{ and } b_3(a,\Sigma)=0.$$
    \end{itemize}
  \item[(7)] If $\boxed{k=3}$, then $d_1(\Sigma)=n-m$, where $m$ is the maximum number of lines $V(\ell_i)\subset \mathbb P^2, i=1,\ldots,t$, counted with multiplicity $m_i$ that pass through the same point.
      \begin{itemize}
    \item[(7.i)] If $1\leq a \leq d_1(\Sigma)$, then $$b_i(a,\Sigma)={{a+2}\choose{a+i-1}}\cdot{{a+i-2}\choose{a-1}}, i=1,2,3.$$
    \item[(7.ii)] If $d_1(\Sigma)+1\leq a \leq d_2(\Sigma)=n-m_1$, then $$b_2(a,\Sigma)=2b_1(a,\Sigma)-a-2 \mbox{ and } b_3(a,\Sigma)=b_1(a,\Sigma)-a-1.$$
    \item[(7.iii)] If $d_2(\Sigma)+1\leq a\leq d_3(\Sigma)=n$, then $$b_2(a,\Sigma)=2b_1(a,\Sigma)-e-2 \mbox{ and } b_3(a,\Sigma)=b_1(a,\Sigma)-e-1.$$
  \end{itemize}
  \item[(8)] If $a\geq d_1(\Sigma)+1$, then ${\rm ht}(I_a(\Sigma))\leq k-1$ and so $R/I_a(\Sigma)$ is Cohen-Macaulay iff $\Sigma$ is $(n-a)$-generic. In this case, $$b_i(a,\Sigma)={{n}\choose{i+a-1}}\cdot{{i+a-2}\choose{a-1}}, i=1,\ldots,n-a+1, \mbox{ and } b_{n-a+2}(a,\Sigma)=0.$$
\end{itemize}
\end{thm}

\section{Betti numbers via deletion-contraction}\label{sec_DelContr}

One of the most common way to compute Betti numbers is to use recursion and {\em mapping cone} technique. Suppose $0\rightarrow M \rightarrow N \rightarrow P\rightarrow 0$ be a short exact sequence of finitely generated $R$-modules, where $R:=\mathbb K[x_1,\ldots,x_k]$. Suppose ${\bf F}_{\bullet}\twoheadrightarrow M$ and ${\bf G}_{\bullet}\twoheadrightarrow N$ are minimal free resolutions\footnote{Working over the polynomial ring $R$ guarantees that these resolutions are finite.} of $M$, and respectively, $N$. Then ${\bf F}_{\bullet-1}\oplus {\bf G}_{\bullet}\twoheadrightarrow P$ is a free resolution (not minimal) of $P$; by convention ${\bf F}_{-1}=0$ and ${\bf F}_{\bullet-1}$ means that the indexing is decreased by 1 (i.e., the $i$-th free module in the free resolution of $P$ is ${\bf F}_{i-1}\oplus {\bf G}_i$). To get a minimal free resolution, one uses ``cancelations'', which means removing direct summands in consecutive terms of the resolution which produce invertible block submatrices of the corresponding maps in the resolution.

Let $\Sigma=(\ell_1,\ldots,\ell_n)\subset R:=\mathbb K[x_1,\ldots,x_k]$ be a collection of linear forms, and let $1\leq a\leq n$ be an integer. Let $\ell\in\Sigma$ be any element, and let $\Sigma':=\Sigma\setminus\{\ell\}$ be the collection of the same linear forms obtained by removing $\ell$ from $\Sigma$. After a change of coordinates, suppose $\ell=x_1$. Let $\bar{R}:=\mathbb K[x_2,\ldots,x_k]$, and let $\bar{\Sigma}:=(\bar{\ell}_2,\ldots,\bar{\ell}_n)\subset \bar{R}$, where $\bar{\ell}_i:=\ell_i(0,x_2,\ldots,x_k)$. In this setup, we also consider $0$ as the ``zero linear form'', and by convention, if $a$ is strictly bigger than the number of nonzero elements of $\bar{\Sigma}$, then $I_a(\bar{\Sigma})=0$, the zero ideal.

\begin{thm}\label{thm_recursion} With the above notations we have that for $i\geq 1$,
$$b_i(a,\Sigma)=b_i(a-1,\Sigma')+b_i(a,\bar{\Sigma})+b_{i-1}(a,\bar{\Sigma}),$$ where $b_0(a,\bar{\Sigma})=0$ by convention, and $b_{k+1}(a,\Sigma) = b_{k+1}(a-1,\Sigma') = b_k(a,\bar{\Sigma})=0$.
\end{thm}
\begin{proof} By \cite[Corollary 2.3]{BuToXi}, we have $I_a(\Sigma):x_1=I_{a-1}(\Sigma')$. Then, we have a short exact sequence of $R$-modules:
$$0\rightarrow \frac{R(-1)}{I_{a-1}(\Sigma')}\rightarrow \frac{R}{I_a(\Sigma)}\rightarrow \frac{R}{\langle I_a(\bar{\Sigma}),x_1\rangle}\rightarrow 0.$$

By \cite[Theorem 2.2]{BuToXi} we have the following linear minimal free resolutions:

(i) For $R(-1)/I_{a-1}(\Sigma')$:

$${\bf F}_{\bullet}\twoheadrightarrow \frac{R(-1)}{I_{a-1}(\Sigma')} \mbox{ with } {\bf F}_i=R^{b_i(a-1,\Sigma')}(-(a+i-1)), {\bf F}_0=R(-1), {\bf F}_j=0, j\geq k+1.$$

(ii) For $R/I_a(\Sigma)$:

$${\bf G}_{\bullet}\twoheadrightarrow \frac{R}{I_{a}(\Sigma)} \mbox{ with } {\bf G}_i=R^{b_i(a,\Sigma)}(-(a+i-1)), {\bf G}_0=R, {\bf G}_j=0, j\geq k+1.$$

(iii) For $R/\langle I_a(\bar{\Sigma}), x_1\rangle$:

$${\bf H}_{\bullet}\twoheadrightarrow \frac{R}{\langle I_a(\bar{\Sigma}), x_1\rangle} \mbox{ with } {\bf H}_i=R^{b_i(a,\bar{\Sigma})+b_{i-1}(a,\bar{\Sigma})}(-(a+i-1)), i\geq 2, {\bf H}_1=R^{b_1(a,\bar{\Sigma})}(-a)\oplus, R(-1),$$ $${\bf H}_0=R, {\bf H}_j=0, j\geq k+1.$$

\medskip

Now, comparing the mapping cone resolution ${\bf F}_{\bullet-1}\oplus {\bf G}_{\bullet}\twoheadrightarrow R/\langle I_a(\bar{\Sigma}), x_1\rangle$ to the minimal resolution ${\bf H}_{\bullet}\twoheadrightarrow R/\langle I_a(\bar{\Sigma}), x_1\rangle$, we get the desired recursive formulas. In more details, comparing $$\cdots\rightarrow R^{b_{i-1}(a-1,\Sigma')}(-(a+i-2))\oplus R^{b_i(a,\Sigma)}(-(a+i-1))\rightarrow R^{b_i(a-1,\Sigma')}(-(a+i-1))\oplus R^{b_{i-1}(a,\Sigma)}(-(a+i-2))\rightarrow\cdots$$ with
$$\cdots\rightarrow R^{b_i(a,\bar{\Sigma})+b_{i-1}(a,\bar{\Sigma})}(-(a+i-1))\rightarrow R^{b_{i-1}(a,\bar{\Sigma})+b_{i-2}(a,\bar{\Sigma})}(-(a+i-2))\rightarrow\cdots$$ the only cancellation possible is the one that gives
$$R^{b_{i-1}(a,\Sigma)-b_{i-1}(a-1,\Sigma')}(-(a+i-2))\cong R^{b_{i-1}(a,\bar{\Sigma})+b_{i-2}(a,\bar{\Sigma})}(-(a+i-2)).$$
\end{proof}

\begin{exm}\label{exm_delcontr} We come back to Example \ref{tutte} and we use Theorem \ref{thm_recursion} to obtain the same Betti numbers we got in Example \ref{exm_betti_tutte}. $\Sigma=(x_1, x_1, x_2, x_3, x_1-x_3, x_2+x_3, x_1+2x_2+5x_3)$. We want to determine $b_i(a,\Sigma),i=1,2,3$, for $a=4,5,6$.

First, consider $\ell=x_1$. Then, $$\Sigma'=(x_1,x_2,x_3,x_1-x_3,x_2+x_3,x_1+2x_2+5x_3) \mbox{ and } \bar{\Sigma}=(0,0,x_2,x_3,x_3,x_2+x_3,2x_2+5x_3),$$ and

\begin{eqnarray}
b_3(a,\Sigma)&=&b_3(a-1,\Sigma')+b_2(a,\bar{\Sigma})\nonumber\\
b_2(a,\Sigma)&=&b_2(a-1,\Sigma')+b_2(a,\bar{\Sigma})+b_1(a,\bar{\Sigma})\nonumber\\
b_1(a,\Sigma)&=&b_1(a-1,\Sigma')+b_1(a,\bar{\Sigma}).\nonumber
\end{eqnarray}

Then, consider $\ell'=x_3$. $$\Sigma''=(x_1,x_2,x_1-x_3,x_2+x_3,x_1+2x_2+5x_3) \mbox{ and } \bar{\Sigma'}=(x_1,x_2,0,x_1,x_2,x_1+2x_2),$$ and

\begin{eqnarray}
b_3(a-1,\Sigma')&=&b_3(a-2,\Sigma'')+b_2(a-1,\bar{\Sigma'})\nonumber\\
b_2(a-1,\Sigma')&=&b_2(a-2,\Sigma'')+b_2(a-1,\bar{\Sigma'})+b_1(a-1,\bar{\Sigma'})\nonumber\\
b_1(a-1,\Sigma')&=&b_1(a-2,\Sigma'')+b_1(a-1,\bar{\Sigma'}).\nonumber
\end{eqnarray}

\medskip

$\bullet$ Observe that $\Sigma''$ is $2$-generic with $|\Sigma''|=5$. If $a=4,5,6$, then $a-2=2,3,4$, and hence $5-(a-2)=3,2,1$. Therefore we have
\begin{itemize}
  \item[(a1)] If $a=4$, because $d(\Sigma'')=5-2=3$, then by Theorem \ref{thm_summary}~(2), $\displaystyle b_i(2,\Sigma'')={{4}\choose{i+1}}\cdot{{i}\choose{1}}, i=1,2,3$.
  \item[(b1)] If $a=5$, then by Theorem \ref{thm_summary}~(8), $\displaystyle b_i(3,\Sigma'')={{5}\choose{i+2}}\cdot{{i+1}\choose{2}}, i=1,2,3$.
  \item[(c1)] If $a=6$, then by Theorem \ref{thm_summary}~(5), $\displaystyle b_1(4,\Sigma'')=5, b_2(4,\Sigma'')=4, \mbox{ and } b_3(4,\Sigma'')=0$.
\end{itemize}

\medskip

$\bullet$ Now we look at $\bar{\Sigma'}=(x_1,x_2,0,x_1,x_2,x_1+2x_2)$. If $a=4,5,6$, then $a-1=3,4,5$. Note that $I_{a-1}(\bar{\Sigma'})=I_{a-1}(x_1,x_1,x_2,x_2,x_1+2x_2)$.
\begin{itemize}
  \item[(a2)] If $a=4$, then by Theorem \ref{thm_summary}~(6.ii), since $e=3$, $b_1(3,\bar{\Sigma'})=4$ and $b_2(3,\bar{\Sigma'})=3$.
  \item[(b2)] If $a=5$, then by Theorem \ref{thm_summary}~(6.ii), since $e=2$, $b_1(4,\bar{\Sigma'})=3$ and $b_2(4,\bar{\Sigma'})=2$.
  \item[(c2)] If $a=6$, then $I_5(\bar{\Sigma'})=\langle x_1^2x_2^2(x_1+2x_2)\rangle$, so $b_1(5,\bar{\Sigma'})=1$ and $b_2(5,\bar{\Sigma'})=0$.
\end{itemize}

\medskip

Withe the two bullets above we get $b_i(a-1,\Sigma')$ as follows:
\begin{itemize}
  \item[(A1)] If $a=4$, then (a1) and (a2) give $$b_1(3,\Sigma')=6+4=10, b_2(3,\Sigma')=8+3+4=15, b_3(3,\Sigma')=3+3=6.$$
  \item[(B1)] If $a=5$, then (b1) and (b2) give $$b_1(4,\Sigma')=10+3=13, b_2(4,\Sigma')=15+2+3=20, b_3(4,\Sigma')=6+2=8.$$
  \item[(C1)] If $a=6$, then (c1) and (c2) give $$b_1(5,\Sigma')=5+1=6, b_2(5,\Sigma')=4+0+1=5, b_3(5,\Sigma')=0+0=0.$$
\end{itemize}

\medskip

$\bullet$ Lastly, we need to analyze the Betti numbers of $\bar{\Sigma}=(0,0,x_2,x_3,x_3,x_2+x_3,2x_2+5x_3)$.
\begin{itemize}
  \item[(A2)] If $a=4$, then $I_4(\bar{\Sigma})=I_4(x_2,x_3,x_3,x_2+x_3,2x_2+5x_3)$. From Theorem \ref{thm_summary}~(6.ii), since $e=3$, we have $b_1(4,\bar{\Sigma})=4$ and $b_2(4,\bar{\Sigma})=3$.
  \item[(B2)] If $a=5$, then $I_5(\bar{\Sigma})=I_5(x_2,x_3,x_3,x_2+x_3,2x_2+5x_3)=\langle x_2x_3^2(x_2+x_3)(2x_2+5x_3)\rangle$. $b_1(5,\bar{\Sigma})=1$ and $b_2(5,\bar{\Sigma})=0$.
  \item[(C2)] If $a=6$, then $I_6(\bar{\Sigma})=0$, and so, $b_1(6,\bar{\Sigma})=b_2(6,\bar{\Sigma})=0$.
\end{itemize}

\medskip

And now we can conclude the calculations, observing that we obtained the same results as in Example \ref{exm_betti_tutte}.

\begin{itemize}
  \item If $\boxed{a=4}$, then (A1) and (A2) give $$b_1(4,\Sigma)=10+4=14, b_2(4,\Sigma)=15+3+4=22, b_3(4,\Sigma)=6+3=9.$$
  \item If $\boxed{a=5}$, then (B1) and (B2) give $$b_1(5,\Sigma)=13+1=14, b_2(5,\Sigma)=20+0+1=21, b_3(5,\Sigma)=8+0=8.$$
  \item If $\boxed{a=6}$, then (C1) and (C2) give $$b_1(6,\Sigma)=6+0=6, b_2(6,\Sigma)=5+0+0=5, b_3(4,\Sigma)=0+0=0.$$
\end{itemize}
\end{exm}

\begin{rem}\label{rem_Veronese} Suppose $\Sigma=(\underbrace{x_1,\ldots,x_1}_{m_1},\ldots,\underbrace{x_k,\ldots,x_k}_{m_k})\subset R=\mathbb K[x_1,\ldots,x_k]$, with $m_1\geq\cdots\geq m_k\geq 1$, and let $1\leq a\leq n:=m_1+\cdots+m_k$. Then, $I_a(\Sigma)$ is generated by all monomials $x_1^{i_1}\cdots x_k^{i_k}$ with $i_1+\cdots+i_k=a$, and $0\leq i_j\leq m_j$, for all $j=1,\ldots,k$. With this, if $a\geq m_1$, then $I_a(\Sigma)$ becomes a {\em Veronese type ideal}, often denoted by $I_{{\bf m},k,a}$, where ${\bf m}:=(m_1,\ldots,m_k)$.

A smart argument in \cite[Theorem 2.2]{AbZa} presents the minimum number of generators for a Veronese type ideal:

$$\beta_1(a,\Sigma)={{a+k-1}\choose{k-1}}+\sum_{J\subseteq [k]}(-1)^{|J|}{{a+k-1-\sum_{i\in J}(m_i+1)}\choose{k-1}}.$$ Furthermore, with this nice description in hand, Theorem \ref{thm_summary}~(3) will give the same results as in Section 3 of \cite{AbZa}, since both approaches use the strength of Herzog-K\"{u}hl equations.

One must mention that this calculation doesn't cover the entire range of $a$. By Theorem \ref{thm_summary}~(0), $d_1(\Sigma)=n-(m_1+\cdots+m_{k-1})=m_k$. If $1\leq a\leq m_k$, $b_i(a,\Sigma)$ is given by Theorem \ref{thm_summary}~(2), so what is left to consider is $b_i(a,\Sigma)$, when $a\geq m_k+1$. In particular, with the result above, one needs a formula for $b_1(a,\Sigma)$, when $m_k+1\leq a\leq m_1-1$. In Corollary \ref{cor_k3generic} below we give the formula for $b_1(a,\Sigma)$, when $k=3$; hence, by Theorem \ref{thm_summary}~(7), the formulas for all Betti numbers.
\end{rem}

We have tried to use our Theorem \ref{thm_recursion} to recover the result presented in the above remark via induction, and we found it to be rather complicated, compared to the proof of that result that uses inclusion-exclusion principle. The advantage of our Theorem \ref{thm_recursion} is that it can be applied to any specific example (like Example \ref{exm_delcontr} above), with no assumption on knowing the underlying combinatorics. Picking any sequence of hyperplanes to delete and contract at, will reduce the standard computations to the case $k=2$, which is straightforwardly presented in Theorem \ref{thm_summary}~(6). This advantage could be put to good use to computing the Betti numbers of $I_a(\Sigma)$ via an algorithm that doesn't appeal to computing the entire minimal free resolution of $I_a(\Sigma)$.

\subsection{The case $k=3$.} \label{sec_k3} We now take a look on how we can apply Theorem \ref{thm_recursion} in the case $k=3$. Indeed, Theorem \ref{thm_summary}~(7) presents all Betti numbers in terms of the first one, and our goal here is to avoid appealing to the combinatorial formula for the first Betti number as expressed in Theorem \ref{thm_summary}~(1), but to produce more direct formulas. The first advantage for working in $k=3$ is that $\bar{\Sigma}$ in Theorem \ref{thm_recursion} is of rank 2, and therefore $b_i(a,\bar{\Sigma})$ can be computed immediately from Theorem \ref{thm_summary}~(6). The other advantage for this case is that we can draw $\Sigma$ and therefore have a visual perspective of the combinatorics involved.

\medskip

Suppose $\Sigma=(\underbrace{\ell_1,\ldots,\ell_1}_{m_1},\ldots,\underbrace{\ell_t,\ldots,\ell_t}_{m_t})\subset R=\mathbb K[x_1,x_2,x_3]$ with $\gcd(\ell_i,\ell_j)=1$ if $i\neq j$, and $m_i\geq 1$, but not necessarily in decreasing order. Suppose we want to delete and contract at the line defined by $\ell_1$, which we can assume, after a change of variables, that equals $x_1$.

Suppose ${\rm ht}(\langle x_1,\ell_{s_u+1},\ldots,\ell_{s_{u+1}}\rangle)=2$, for $u=0,\ldots,v-1$, with $s_0=1$ and $s_v=t$; this condition is capturing how the remaining $n-m_1$ lines (counting their multiplicities) intersect the line $V(x_1)$. Namely, for $u=0,\ldots,v-1$, $\{P_{u+1}\}=V(x_1)\cap V(\ell_{s_u+1})\cap\cdots\cap V(\ell_{s_{u+1}})$ are all the singular points on $V(x_1)$.

This implies that in $\bar{\Sigma}$, $\ell'_{u+1}:\equiv \ell_{s_u+1}\equiv\cdots\equiv\ell_{s_{u+1}} \mbox{ mod } x_1$ will occur with multiplicity $n_{u+1}:=m_{s_u+1}+\cdots+m_{s_{u+1}}$, and $\gcd(\ell'_{u+1},\ell'_{u'+1})=1$ for any distinct $u, u'\in\{0,\ldots,v-1\}$.\footnote{It is worth noting that $m$ in Theorem \ref{thm_summary}~(7), is greater than or equal to $m_1+\max\{n_1,\ldots,n_v\}$.} Therefore,

$$\bar{\Sigma}=(\underbrace{\ell'_1,\ldots,\ell'_1}_{n_1},\ldots,\underbrace{\ell'_v,\ldots,\ell'_v}_{n_v})\subset\mathbb K[x_1,x_2],$$ and apply Theorem \ref{thm_summary}~(6), after reordering the $\ell'_j$ such that the multiplicities will occur in decreasing order, to find the Betti numbers of $I_a(\bar{\Sigma})$. If $a>n-m_1=n_1+\cdots+n_v$, then $I_a(\bar{\Sigma})=0$. On the other hand, if $a\leq n_1+\cdots+n_v-\max\{n_1,\ldots,n_v\}$, then we are in Case (6.i) of Theorem \ref{thm_summary}, and we know the Betti numbers of the contraction once again.

We then use these in Theorem \ref{thm_recursion}, with $\Sigma'=(\underbrace{\ell_1,\ldots,\ell_1}_{m_1-1},\underbrace{\ell_2,\ldots,\ell_2}_{m_2},\ldots,\underbrace{\ell_t,\ldots,\ell_t}_{m_t})$. Maybe it is convenient to delete and contract at $\ell_1$ again, until we get rid off it with all its multiplicity. Therefore we get that for $i=1,2,3$,

$$b_i(a,\Sigma)=b_i(a-m_1,\tilde{\Sigma})+\sum_{j=0}^{m_1-1}[b_i(a-j,\bar{\Sigma})+b_{i-1}(a-j,\bar{\Sigma})],$$ where $\tilde{\Sigma}=(\underbrace{\ell_2,\ldots,\ell_2}_{m_2},\ldots,\underbrace{\ell_t,\ldots,\ell_t}_{m_t})$.

\medskip

Now we are filling the gaps mentioned in Remark \ref{rem_Veronese} using the above strategy. By Theorem \ref{thm_summary}~(7), it is enough to consider $i=1$ and $m_3+1\leq a\leq m_2+m_3$, if we assume $m_1\geq m_2\geq m_3\geq 1$. By Remark \ref{rem_Veronese}, we may assume $a\leq m_1-1$ and so, $m_1\geq m_3+2$.

\begin{cor}\label{cor_k3generic} $\Sigma=(\underbrace{x_1,\ldots,x_1}_{m_1},\underbrace{x_2,\ldots,x_2}_{m_2},\underbrace{x_3,\ldots,x_3}_{m_3})\subset R=\mathbb K[x_1,x_2,x_3]$, with $m_1\geq m_2\geq m_3\geq 1$, $m_1\geq m_3+2$, and let $n:=m_1+m_2+m_3$. Let $a$ be an integer with $m_3+1\leq a\leq m_2+m_3$ and $a\leq m_1-1$. Then,

\begin{itemize}
  \item[(i)] If $a\leq m_2$, then $$b_1(a,\Sigma)=\sum_{j=0}^{m_3}(a-j+1)=(m_3+1)(a+1)-{{m_3+1}\choose{2}}.$$
  \item[(ii)] If $a\geq m_2+1$, then
  \begin{eqnarray}
  b_1(a,\Sigma)&=&\sum_{j=a-m_2}^{m_3}(a-j+1)+\sum_{j=0}^{a-m_2-1}(m_2+1)\nonumber\\
  &=&(a+1)(n-m_1-a+1)+(a-m_2)(m_2+1)+{{a-m_2}\choose{2}}-{{m_3+1}\choose{2}}.\nonumber
  \end{eqnarray}
\end{itemize}
\end{cor}
\begin{proof} We use the above strategy by deleting and contracting at the line $V(x_3)$. Since $\bar{\Sigma}=\tilde{\Sigma}=(\underbrace{x_1,\ldots,x_1}_{m_1},\underbrace{x_2,\ldots,x_2}_{m_2})$, we have

$$b_1(a,\Sigma)=\sum_{j=0}^{m_3}b_1(a-j,\bar{\Sigma}).$$

Next we apply Theorem \ref{thm_summary}~(6), hence we need to see how $a-j$ compares to $d_1(\bar{\Sigma})=m_1+m_2-m_1=m_2$.

\medskip

$\bullet$ \underline{Case $a\leq m_2$}. In this case, $a-j\leq m_2$, for all $j=0,\ldots,m_3$, hence $b_1(a-j,\bar{\Sigma})=a-j+1$. This leads to formula (i).

\medskip

$\bullet$ \underline{Case $a\geq m_2+1$}. If $a-m_2\leq j\leq m_3$, then $m_2\geq a-j\geq a-m_3$, and so, for this range of $j$, we have $b_1(a-j,\bar{\Sigma})=a-j+1$.

Suppose $0\leq j\leq a-m_2-1$. Then $e_1=\max\{m_1+a-j-(m_1+m_2),0\}=a-j-m_2$, as $a-j\geq m_2+1$. On the other hand, $e_2=\max\{m_2+a-j-(m_1+m_2),0\}=0$, as $a\leq m_1-1$. So $e=\max\{a-j-(e_1+e_2),0\}=m_2$, leading that $b_1(a-j,\bar{\Sigma})=m_2+1$, for this range of $j$. This leads to formula (ii).
\end{proof}

\medskip

Another situation when we can obtain compact formulas for the Betti numbers is the following. Suppose $\Sigma=(\ell_1,\ldots,\ell_n)\subset R=\mathbb K[x_1,x_2,x_3], n\geq 4$ defines a line arrangement $\mathcal A$ of rank 3 in $\mathbb P^2$; not that the multiplicity of each line is one. Let $m$ be the invariant presented in Theorem \ref{thm_summary}~(7), i.e., the maximum number of concurent lines of $\mathcal A$.\footnote{We can assume $m\geq 3$, since if $m=2$, then $\A$ becomes 2-generic (i.e., generic line arrangement) and the Betti numbers are explicitly presented at the end of Section \ref{CM}.} After a change of variables, suppose $\ell_1=x_1$ and that $P_1,\ldots,P_t$ are all the singularities of $\mathcal A$ each having $m$ lines of $\A$ (including, possibly $V(\ell_1)$) passing through it. In these conditions, $d:=d_1(\Sigma)=n-m$.

\begin{cor}\label{cor_dplus1} With the above notations, assume furthermore that $\{P_1,\ldots,P_t\}\subset V(\ell_1)$. Then,
$$b_1(n-m+1,\Sigma)={{n-m+3}\choose{2}}-t.$$
\end{cor}
\begin{proof} Consider $\Sigma'=\Sigma\setminus\{\ell_1\}$. Then $|\Sigma'|=n-1=:n'$, and the maximum number of concurent lines of the arrangement defined by $\Sigma'$ is $m'=m-1$. Also, the rank of $\Sigma'$ remains 3. With these, $d_1(\Sigma')=n'-m'=n-m=d$. Then, by Theorem \ref{thm_summary}~(2), $\displaystyle b_1(d,\Sigma')={{n-m+2}\choose{2}}$.

Suppose $Q_{t+1},\ldots,Q_s$ are the remaining singularities of $\mathcal A$ placed on the line $V(\ell_1)$. Each point $Q_j$ has $m_j$ lines, including $V(\ell_1)$, passing through it. By assumptions, $m_j<m$. Then the contraction at $\ell_1$ will produce $\bar{\Sigma}$ which consists of $s$ nonproportional forms in $\mathbb K[x_2,x_3]$ with multiplicities $n_1=\cdots=n_t=m-1$ and $n_j=m_j-1, j=t+1,\ldots,s$. Also, $|\bar{\Sigma}|=n-1$.

We want to calculate $b_1(n-m+1,\bar{\Sigma})$. For $i=1,\ldots,t$, $e_i=\max\{(m-1)+(n-m+1)-(n-1),0\}=1$, and for $j=t+1,\ldots,s$, $e_j=\max\{(m_j-1)+(n-m+1)-(n-1),0\}=0$. So $e=\max\{(n-m+1)-t,0\}=n-m+1-t$, giving that $b_1(n-m+1,\bar{\Sigma})=n-m+2-t$, by Theorem \ref{thm_summary}~(6).

From Theorem \ref{thm_recursion} we obtain

$$b_1(n-m+1,\Sigma)=b_1(n-m,\Sigma')+b_1(n-m+1,\bar{\Sigma})={{n-m+2}\choose{2}}+n-m+2-t={{n-m+3}\choose{2}}-t.$$
\end{proof}

\begin{exm}\label{exm_dplus1} Let $\Sigma=(x_1,x_1-x_3,x_1-2x_3,x_1-3x_3,x_2,x_1-x_2,x_1-2x_2,x_1+x_2-2x_3)\subset R=\mathbb K[x_1,x_2,x_3]$. The (projective) picture is the following:

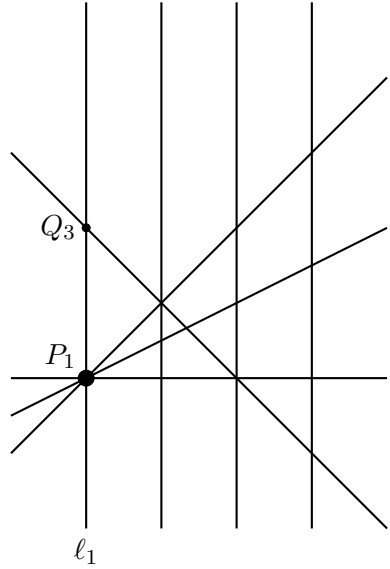
\begin{center}
\begin{tikzpicture}
\draw[thick] (0,-2) node[below]{$\ell_1$} -- (0,5);
\draw[thick] (1,-2) -- (1,5);
\draw[thick] (2,-2) -- (2,5);
\draw[thick] (3,-2) -- (3,5);
\draw[thick] (-1,0) -- (4,0);
\draw[thick] (-1,-1) -- (4,4);
\draw[thick] (-1,-0.5) -- (4,2);
\draw[thick] (-1,3) -- (4,-2);
\filldraw[black] (0,0) node[above left]{$P_1$} circle (3pt);
\filldraw[black] (1.5,6) node[left]{$P_2$} circle (3pt);
\filldraw[black] (0,2) node[left]{$Q_3$} circle (1.5pt);
\end{tikzpicture}
\end{center}

We have $n=8$, $m=4$, $t=2$. Then $\displaystyle b_1(4+1,\Sigma)={{4+3}\choose{2}}-2=19$, which matches computations performed with \cite{GrSt}. Also, note that we can replace the line $V(x+y-2z)$ with any line that doesn't pass through $P_1$ nor $P_2$; the result would be the same.
\end{exm}

\section{Combinatorial formulas for Betti numbers}\label{sec_comb}

Let $\Sigma=(\ell_1,\ldots,\ell_n)\subset R:=\mathbb K[x_1,\ldots,x_k]$ be a collection of linear forms, and let $1\leq a\leq n-1$ be an integer. Suppose $\langle \ell_1,\ldots,\ell_n\rangle=\langle x_1,\ldots,x_k\rangle$. Suppose we fixed the defining linear forms $\ell_i$. Knowledge of the Hilbert function also can help calculating Betti numbers in a recursive way.

\begin{lem}\label{lem_HF} For $i=1,\ldots,k$ one has

$$b_i(a,\Sigma)=\sum_{j=1}^{i-1}(-1)^{j-1}{{k+j-1}\choose{j}}b_{i-j}(a,\Sigma)+(-1)^{i-1}\underbrace{\dim_{\mathbb K}(I_a(\Sigma))_{a+i-1}}_{{\rm HF}(I_a(\Sigma), a+i-1)}.$$

\end{lem}
\begin{proof}
We have the linear free resolution

$$0\longrightarrow R^{b_{k}(a,\Sigma)}(-(a+k-1))\longrightarrow \cdots \longrightarrow R^{b_2(a,\Sigma)}(-(a+1))\longrightarrow R^{b_1(a,\Sigma)}(-a)\longrightarrow I_a(\Sigma)\rightarrow 0.$$ Applying the Hilbert function in degree $a+i-1$, and since $\displaystyle{\rm HF}(R,j)=\dim_{\mathbb K}R_j={{k+j-1}\choose{j}}$, we obtain the conclusion.
\end{proof}

\subsection{An exact complex.} To proceed further, we make the following notations, some more standard than the others: $[n]:=\{1,\ldots,n\}$, and if $J=\{j_1,\ldots,j_s\}\subseteq [n]$, we denote $\displaystyle F_J:=\prod_{i\in J}\ell_i=\ell_{j_1}\cdots\ell_{j_s}$. So $$I_a(\Sigma)=\langle F_J|J\subseteq [n], |J|=a\rangle.$$ Also, if $J=\{j_1,\ldots,j_s\}\subseteq [n]$, $r(J):=\dim_{\mathbb K}{\rm Span}\{\ell_{j_1},\ldots,\ell_{j_s}\}$.

Consider the following surjective $R$-modules homomorphism:

$$\Phi_{a,\Sigma}:\bigoplus_{1\leq j_1<\cdots<j_{n-a}\leq n}\frac{R(-a)}{\langle \ell_{j_1},\ldots,\ell_{j_{n-a}}\rangle}\longrightarrow \frac{I_a(\Sigma)}{I_{a+1}(\Sigma)},$$ given by $$\Phi_{a,\Sigma}(\ldots,\hat{h}_{j_1,\ldots,j_{n-a}},\ldots)=(\sum h_{j_1,\ldots,j_{n-a}}F_{[n]\setminus\{j_1,\ldots,j_{n-a}\}}) \mbox{ mod } I_{a+1}(\Sigma).$$ By $\hat{f}$ we understand the residue class of an $f\in R$ modulo the ideal $\langle \ell_{j_1},\ldots,\ell_{j_{n-a}}\rangle$ corresponding to the location $(j_1,\ldots,j_{n-a})$ of $f$.

$\Phi_{a,\Sigma}$ is well-defined: if $h_{j_1,\ldots,j_{n-a}}=g_1\ell_{j_1}+\cdots+g_{n-a}\ell_{j_{n-a}}$, then
$$h_{j_1,\ldots,j_{n-a}}F_{[n]\setminus\{j_1,\ldots,j_{n-a}\}}=g_1F_{[n]\setminus\{j_2,\ldots,j_{n-a}\}}+\cdots+ g_{n-a}F_{[n]\setminus\{j_1,\ldots,j_{n-a-1}\}}\in I_{a+1}(\Sigma).$$

\medskip

Let $1\leq s\leq n-a+1$ be an integer, and suppose $\{j_1,\ldots,j_s\}$ is a circuit (i.e., minimal dependent set) of length $s$ among the linear forms $\ell_{j_1},\ldots,\ell_{j_s}$. Since we have the $\ell_i$'s fixed, this comes with a unique (up to multiplication by a nonzero scalar) linear dependency
$$D_{j_1,\ldots,j_s}: c_{j_1}\ell_{j_1}+\cdots+c_{j_s}\ell_{j_s}=0, c_{j_i}\in\mathbb K\setminus\{0\}, i=1,\ldots,s.$$

Multiplying this by $F_{[n]\setminus\{j_1,\ldots,j_s\}}$, one obtains the relation.

$$c_{j_1}F_{[n]\setminus\{j_2,\ldots,j_s\}}+\cdots+c_{j_s}F_{[n]\setminus\{j_1,\ldots,j_{s-1}\}}=0.$$ Since $|[n]\setminus\{j_1,\ldots,j_s\}|=n-s\geq a-1$, then there are $n-s$ linear forms of $\Sigma$ that are common factors of every $G_i:=F_{([n]\setminus\{j_1,\ldots,j_s\})\cup\{j_i\}}$, $i=1,\ldots,s$. But each $G_i\in I_{n-s+1}(\Sigma)$, so dividing every $G_i$ by the same product of $n-s+1-a\geq 0$ common linear forms, we get a relation with constant nonzero coefficients among the generators of $I_a(\Sigma)$.

So the dependency $D_{j_1,\ldots,j_s}$ will produce $\displaystyle {{n-s}\choose{n-s-a+1}}$ distinct vectors $${\bf c}:=(0,\ldots,0,c_{j_1},0,\ldots,0,c_{j_i},0,\ldots,0,c_{j_{i+1}},0,\ldots,0,c_{j_s},0,\ldots,0)\in \ker\Phi_{a,\Sigma}.$$ The locations of the nonzero entries $c_{j_i}$ depend on the simplification by the $n-s+1-a$ common linear forms used. So two nonzero entries in two distinct such vectors will never be on the same location in their vectors. This implies that all these vectors are linearly independent. So, each circuit $\{j_1,\ldots,j_s\}$ of $\Sigma$ produces a $\mathbb K$-linear subspace $\mathcal W_{\{j_1,\ldots,j_s\}}(a,\Sigma)\subset \ker\Phi_{a,\Sigma}$ of dimension $\displaystyle {{n-s}\choose{n-s-a+1}}$.

\medskip

Let $$\mathcal W(a,\Sigma):= \sum_{s=1}^{n-a+1}\sum_{J \mbox{ circuit of length }s}\mathcal W_{\{j_1,\ldots,j_s\}}(a,\Sigma),$$ and suppose $$N_{a,\Sigma}:=\dim_{\mathbb K}\mathcal W(a,\Sigma).$$ Since multiplying any vector of $\mathcal W(a,\Sigma)$ by any homogeneous polynomial, one obtains another element of $\ker\Phi_{a,\Sigma}$, if ${\bf c}_1,\ldots, {\bf c}_{N_{a,\Sigma}}$ is a basis for $\mathcal W(a,\Sigma)$ that we place as the rows of an $\displaystyle N_{a,\Sigma}\times {{n}\choose{n-a}}$ matrix $\mathcal M$, we obtain an $R$-homomorphism
$$\Psi_{a,\Sigma}: R^{N_{a,\Sigma}}\stackrel{\cdot \mathcal M}\longrightarrow \bigoplus_{1\leq j_1<\cdots<j_{n-a}\leq n}\frac{R(-a)}{\langle \ell_{j_1},\ldots,\ell_{j_{n-a}}\rangle},$$ whose image we denote with $\mathcal Y(a,\Sigma)$ is included in $\ker\Phi_{a,\Sigma}$.

So we get the complex of $R$-modules

$$(*)\,\, 0\rightarrow \ker\Psi_{a,\Sigma}(-a)\rightarrow R^{N_{a,\Sigma}}(-a)\stackrel{\Psi_{a,\Sigma}}\longrightarrow \bigoplus_{1\leq j_1<\cdots<j_{n-a}\leq n}\frac{R(-a)}{\langle \ell_{j_1},\ldots,\ell_{j_{n-a}}\rangle}\stackrel{\Phi_{a,\Sigma}}\longrightarrow \frac{I_a(\Sigma)}{I_{a+1}(\Sigma)}\rightarrow 0.$$

\begin{rem}\label{rem_kernelPsi} (a) It is interesting to see what are the possible columns of $\mathcal M$ that are all zero, because these will not impose any conditions on the dimension of $\ker(\Psi_{a,\Sigma})$ in any degree $i\geq 0$. The claim is that the column of $\mathcal M$ in the position $(j_1,\ldots,j_{n-a})$ is zero if and only if for any $u\in [n]\setminus\{j_1,\ldots,j_{n-a}\}$, we have $\ell_u\notin\langle \ell_{j_1},\ldots,\ell_{j_{n-a}}\rangle$.

To see this, we have $\ell_u\in \langle \ell_{j_1},\ldots,\ell_{j_{n-a}}\rangle$ iff without loss of generality, $\ell_u=c_{j_1}\ell_{j_1}+\cdots+c_{j_{s-1}}\ell_{j_{s-1}}$, for some $1\leq s-1\leq n-a$ and with $\{j_1,\ldots,j_{s-1},u\}$ being a circuit of length $s\leq n-a+1$. Then we obtain the relation $$c_{j_1}F_{[n]\setminus\{j_2,\ldots,u\}}+\cdots+(-1)F_{[n]\setminus\{j_1,\ldots,j_{s-1}\}}=0,$$ that will be divided by the product $\ell_{j_s}\cdots\ell_{j_{n-a}}$. This lead to $-1$ being placed on the $(j_1,\ldots,j_{n-a})$-th column of $\mathcal M$. The converse is immediate, by reversing what we just have done.

If we do a change of variables and if $r(\{j_1,\ldots,j_{n-a}\})=\delta$, then we can assume $\Sigma_1:=(\ell_{j_1},\ldots,\ell_{j_{n-a}})\subset\mathbb K[x_1,\ldots, x_{\delta}]$ and, for any $i\in[n]\setminus\{j_1,\ldots,j_{n-a}\}$, we have $\ell_i\in\mathbb K[x_{\delta+1},\ldots,x_k]$. In other words, $\Sigma=\Sigma_1\oplus (\Sigma\setminus\Sigma_1)$.

\medskip

(b) If $r(\{j_1,\ldots,j_{n-a}\})=k$, then the $(j_1,\ldots,j_{n-a})$ column on $\mathcal M$ will not impose any condition on the dimension of $\ker(\Psi_{a,\Sigma})$ in any degree $i\geq 1$, since any linear combination with coefficients the entries of this column of homogeneous polynomials of degree $i$ will produce a homogeneous polynomial of degree $i$ that is always an element of $\langle \ell_{j_1},\ldots,\ell_{j_{n-a}}\rangle=\frak m$.
\end{rem}

\medskip

\begin{thm}\label{thm_kernel} With the notations above, if $1\leq a\leq n-1$, then:
\begin{itemize}
  \item[(i)] $\dim_{\mathbb K}(\ker\Phi_{a,\Sigma})_0={{n}\choose{n-a}}-b_1(a,\Sigma)$.
  \item[(ii)] $\ker\Phi_{a,\Sigma}=\mathcal Y(a,\Sigma)$, and so, the complex $(*)$ is exact.
  \item[(iii)] $N_{a,\Sigma}={{n}\choose{n-a}}-b_1(a,\Sigma)$.
  \item[(iv)] $(\ker\Psi_{a,\Sigma})_0=0$, and for any $i\geq 1$, we have
  \begin{eqnarray}
  {\rm HF}(I_a(\Sigma), a+i)-{\rm HF}(I_{a+1}(\Sigma),a+i)&=&\sum_{1\leq j_1<\cdots<j_{n-a}\leq n} {{k-r(\{j_1,\ldots,j_{n-a}\})-1+i}\choose{i}}\nonumber\\
  &-&N_{a,\Sigma}{{k-1+i}\choose{i}}+\dim_{\mathbb K}(\ker\Psi_{a,\Sigma})_i.\nonumber
  \end{eqnarray}
\end{itemize}
\end{thm}
\begin{proof} (i) Looking at the graded $R$-homomorphism $\Phi_{a,\Sigma}$ in degree $a$, since $(I_{a+1}(\Sigma))_a=0$ and $R_0=\mathbb K$, we have the short exact sequence of $\mathbb K$-vector spaces:

$$0\longrightarrow (\ker\Phi_{a,\Sigma})_0\longrightarrow \bigoplus_{1\leq j_1<\cdots<j_{n-a}\leq n}\mathbb K\longrightarrow (I_a(\Sigma))_a\longrightarrow 0.$$ With $b_1(a,\Sigma)=\dim_{\mathbb K}(I_a(\Sigma))_a$, we conclude the proof.

\medskip

(ii) We prove the result by induction on $n:=|\Sigma|\geq 2$.

\medskip

If $\boxed{n=2}$, then $a=1$, $\Sigma=(\ell_1,\ell_2)$, and $$\Phi_{1,\Sigma}: \frac{R(-1)}{\langle \ell_1\rangle}\oplus \frac{R(-1)}{\langle \ell_1\rangle}\longrightarrow \frac{\langle \ell_1,\ell_2\rangle}{\langle \ell_1\ell_2\rangle}.$$ Then, $(\hat{h}_1,\hat{h}_2)\in\ker\Phi_{1,\Sigma}$ iff $h_1\ell_2+h_2\ell_1=h\ell_1\ell_2$, for some $h\in R$.

If $\gcd(\ell_1,\ell_2)=1$, then $\ell_1|h_1$ and $\ell_2|h_2$, and hence $\ker\Phi_{1,\Sigma}=0$. At the same time, there is no linear dependency between $\ell_1$ and $\ell_2$, hence $\mathcal W(1,\Sigma)=0$ and therefore $\mathcal Y(1,\Sigma)=0$.

If $\ell_2=c\ell_1$, for some $c\in\mathbb K\setminus\{0\}$, then $ch_1+h_2=ch\ell_1$. Then $$(\hat{h}_1,\hat{h}_2)=(\hat{h}_1,-c\hat{h}_1)=h_1(\hat{1},-\hat{c})\in\mathcal Y(1,\Sigma),$$ since this vector comes from the dependency $1\cdot\ell_2+(-c)\ell_1=0$. So, $\ker\Phi_{1,\Sigma}\subseteq\mathcal Y(1,\Sigma)$, hence equality.

\medskip

Suppose $\boxed{n>2}$. Let ${\bf \hat{h}}:=(\ldots,\hat{h}_{j_1,\ldots,j_{n-a}},\ldots)\in\ker\Phi_{a,\Sigma}$, where $h_{j_1,\ldots,j_{n-a}}\in R$ of degree $d$; we will say that the vector ${\bf \hat{h}}$ has degree $d$. Then
$${\bf P}:=\sum_{1\leq j_1<\cdots<j_{n-a}\leq n} h_{j_1,\ldots,j_{n-a}}F_{[n]\setminus\{j_1,\ldots,j_{n-a}\}}\in I_{a+1}(\Sigma).$$

Let $\ell_n\in\Sigma$ and suppose we are ordering the terms of ${\bf \hat{h}}$ such that the first $u:=\displaystyle {{n-1}\choose{n-a}}$ correspond to the tuples $(j_1,\ldots,j_{n-a})$ with $1\leq j_1<\cdots<j_{n-a}\leq n-1$, and the next (and last) $v:=\displaystyle {{n-1}\choose{n-a-1}}$ correspond to tuples $(j_1,\ldots,j_{n-a-1},n)$ with $1\leq j_1<\cdots<j_{n-a-1}\leq n-1$. With this, decompose the vector ${\bf \hat{h}}$ as $${\bf \hat{h}}=({\bf \hat{h'}},\underbrace{0,\ldots,0}_v)+(\underbrace{0,\ldots,0}_u,{\bf \hat{h''}}).$$

Also, we can decompose ${\bf P}={\bf A}+{\bf B}$, where

$${\bf A}:=\sum_{1\leq j_1<\cdots<j_{n-a}\leq n-1} h_{j_1,\ldots,j_{n-a}}\ell_nF_{[n-1]\setminus\{j_1,\ldots,j_{n-a}\}}$$ and $${\bf B}:=\sum_{1\leq j_1<\cdots<j_{n-a-1}\leq n-1} h_{j_1,\ldots,j_{n-a-1},n}F_{[n-1]\setminus\{j_1,\ldots,j_{n-a-1}\}}.$$

Let $\Sigma'=\Sigma\setminus\{\ell_n\}$. Let $n':=|\Sigma'|=n-1$. Then we have
\begin{eqnarray}
0=\Phi_{a,\Sigma}({\bf \hat{h}})&=&\Phi_{a,\Sigma}({\bf \hat{h'}},0,\ldots,0)+\Phi_{a,\Sigma}(0,\ldots,0,{\bf \hat{h''}})\nonumber\\
&=&\ell_n\Phi_{a-1,\Sigma'}({\bf \hat{h'}})+\Phi_{a,\Sigma'}({\bf \hat{h''}})\nonumber\\
&=&\Phi_{a-1,\Sigma'}(\ell_n{\bf \hat{h'}})+\Phi_{a,\Sigma'}({\bf \hat{h''}}).\nonumber
\end{eqnarray}

It is clear that ${\bf B}\in I_a(\Sigma')$.

Since ${\bf P}\in I_{a+1}(\Sigma)=\ell_n I_a(\Sigma')+I_{a+1}(\Sigma')\subset I_a(\Sigma')$, we have that ${\bf A}\in I_a(\Sigma')$. So, $\ell_n{\bf \hat{h'}}\in \ker\Phi_{a-1,\Sigma'}$, which is a vector of degree $d+1$. By induction, $\ell_n{\bf \hat{h'}}\in \mathcal Y(a-1,\Sigma')$, and so $$\ell_n{\bf \hat{h'}}=g_1{\bf \hat{c}_1}+\cdots+g_m{\bf \hat{c}_m},$$ for some $g_1,\ldots,g_m\in R_{d+1}$ and ${\bf \hat{c}_1},\ldots,{\bf \hat{c}_m}\in \mathcal W(a-1,\Sigma')$.

Without loss of generality, suppose $\ell_n=x_k$. Then, for $j=1,\ldots,m$, $g_j=\ell_ng'_j+g''_j$, with $g'_j\in R_d$ and $g''_j\in \mathbb K[x_1,\ldots,x_{k-1}]_{d+1}$. Then $$\ell_n({\bf \hat{h'}}-g'_1{\bf \hat{c}_1}-\cdots-g'_m{\bf \hat{c}_m})=g''_1{\bf \hat{c}_1}+\cdots+g''_m{\bf \hat{c}_m}.$$ Since $\ell_n=x_k$ doesn't appear in the righthand-side, then the righthand-side must be zero, and therefore $${\bf\hat{h'}}\in \mathcal Y(a-1,\Sigma').$$ Since $\mathcal W(a-1,\Sigma')$ is built from circuits of length $s\leq n'-(a-1)+1=n-a+1$, then we have $$({\bf \hat{h'}},0,\ldots,0)\in\mathcal Y(a,\Sigma).$$

\medskip

On the other hand, $$\Phi_{a,\Sigma}(0,\ldots,0,{\bf \hat{h''}})=\Phi_{a,\Sigma}({\bf \hat{h}})-\Phi_{a,\Sigma}(\ell_n{\bf \hat{h'}},0,\ldots,0)=0,$$ as $\mathcal Y(a,\Sigma)\subset \ker\Phi_{a,\Sigma}$. So $${\bf \hat{h''}}\in\ker\Phi_{a,\Sigma'}.$$

If $a\leq n'-1$, by induction we have $\ker\Phi_{a,\Sigma'}=\mathcal Y(a,\Sigma')$. But $\mathcal W(a,\Sigma')$ is built from circuits of length $s\leq n'-a+1=n-a<n-a+1$, so $\mathcal W(a,\Sigma')\subseteq \mathcal W(a,\Sigma)$, and so $$(0,\ldots,0,{\bf \hat{h''}})\in \mathcal Y(a,\Sigma).$$ Everything put together gives $${\bf \hat{h}}=({\bf \hat{h'}},0,\ldots,0)+(0,\ldots,0,{\bf \hat{h''}})\in \mathcal Y(a,\Sigma),$$ hence $\ker\Phi_{a,\Sigma}=\mathcal Y(a,\Sigma)$.

\medskip

If $a=n'=n-1$, then the we have

$$\Phi_{n-1,\Sigma}:\frac{R(-(n-1))}{\langle \ell_1\rangle}\oplus\cdots\oplus\frac{R(-(n-1))}{\langle \ell_n\rangle}\rightarrow \frac{I_{n-1}(\Sigma)}{I_n(\Sigma)}$$ given by $$\Phi_{n-1,\Sigma}(\hat{h}_1,\ldots,\hat{h}_n)=h_1\ell_2\cdots\ell_n+\cdots+h_n\ell_1\cdots\ell_{n-1}.$$

$(\hat{h}_1,\ldots,\hat{h}_n)\in\ker\Phi_{n-1,\Sigma}$ iff $h_1\ell_2\cdots\ell_n+\cdots+h_n\ell_1\cdots\ell_{n-1}=h\ell_1\cdots\ell_n$, for some $h\in R$.

It is clear that if $\gcd(\ell_i,\ell_j)=1$ for any $i\neq j$, then $\ell_i|h_i$ for all $i=1,\ldots,n$, and so $\ker\Phi_{n-1,\Sigma}=0$. This matches with the fact that $\mathcal Y(n-1,\Sigma)$ is built from circuits of length $s\leq n-(n-1)+1=2$, and there aren't any such circuits, so $\mathcal Y(n-1,\Sigma)=0$.

Suppose $\Sigma=(\underbrace{\ell_1,\ldots,\ell_1}_{m_1},\ldots,\underbrace{\ell_t,\ldots,\ell_t}_{m_t})$, with $m_1\geq 2$, $m_1\geq m_2\geq\cdots\geq m_t\geq 1$, $t\geq 1$, and $\gcd(\ell_i,\ell_j)=1$, $i\neq j$. Then,

$$\underbrace{(h_1+\cdots+h_{m_1})}_{g_1}\ell_1^{m_1-1}\ell_2^{m_2}\cdots\ell_t^{m_t}+\cdots+ \underbrace{(h_{m_{t-1}+1}+\cdots+h_{m_t})}_{g_t}\ell_1^{m_1-1}\ell_2^{m_2}\cdots\ell_t^{m_t}=h\ell_1^{m_1}\cdots\ell_t^{m_t},$$ leading to $$\ell_1|g_1,\ldots,\ell_t|g_t.$$

With this, the vector $(\hat{h}_1,\hat{h}_2,\ldots,\hat{h}_{m_1},0,\ldots,0)$ can be written as the sum of the vectors
\begin{eqnarray}
(\hat{h}_1,-\hat{h}_1,0,\ldots,0)\nonumber\\
(0,\hat{h}_2+\hat{h}_1,-(\hat{h}_2+\hat{h}_1),0,\ldots,0)\nonumber\\
\vdots\nonumber\\
(0,\ldots,0,\hat{h}_{m_1-1}+\cdots+\hat{h}_1,-(\hat{h}_{m_1-1}+\cdots+\hat{h}_1),0\ldots,0),\nonumber
\end{eqnarray} since $\hat{h}_{m_1}=-(\hat{h}_{m_1-1}+\cdots+\hat{h}_1)$ as $\ell_1|g_1$.

But each such vector obviously belongs to $\mathcal Y(n-1,\Sigma)$ since it corresponds to the dependency $\ell_1-\ell_1=0$ that gives each circuit of length 2 $\{i,i+1\}$, $i=1,\ldots,m_1-1$. Doing the same for the other ``blocks'' of $\ell_2$'s all the way to $\ell_t$'s, we obtain that $(\hat{h}_1,\ldots,\hat{h}_n)\in \mathcal Y(n-1,\Sigma)$, and so, $\ker\Phi_{n-1,\Sigma}=\mathcal Y(n-1,\Sigma)$.

\medskip

(iii) In degree $0$, the graded $R$-homomorphism $\Psi_{a,\Sigma}$ becomes the surjective $\mathbb K$-linear map $$\mathbb K^{N_{a,\Sigma}}\stackrel{\cdot\mathcal M}\longrightarrow \mathcal Y(a,\Sigma)_0=\mathcal W(a,\Sigma),$$ which is an isomorphism of $\mathbb K$-vector spaces, as the rows of $\mathcal M$ are $\mathbb K$-linearly independent vectors. Now, parts (i) and (ii) prove the claim.

\medskip

(iv) This is immediate from (iii), and from looking at the exact complex $(*)$ in degree $a+i$ and apply the previous results.
\end{proof}

\begin{exm}\label{exm_kernelPsi} This is an example where we will see how Theorem \ref{thm_kernel} above can be applied to obtain the Betti numbers. Let $\Sigma=(x_1,x_1,x_1,x_2,x_2)\subset R=\mathbb K[x_1,x_2]$, and let $a=3$. We have $m_1=3$ and $m_2=2$, and $n-m_1=2\leq 3=a$. By Theorem \ref{thm_summary} (6.ii), since $e=2$, we have $b_1(3,\Sigma)=3$, $b_2(3,\Sigma)=2$, and $b_3(3,\Sigma)=0$; these are the numbers we want to obtain, but with this new result.

We have the dependencies corresponding the circuits $\{1,2\}, \{1,3\}, \{2,3\}, \{4,5\}$, which are all the circuits of $\Sigma$:

$$\ell_1-\ell_2=0, \ell_1-\ell_3=0, \ell_2-\ell_3=0, \ell_4-\ell_5=0.$$ Multiplying by the corresponding $F_{[5]\setminus\{i,j\}}$ and simplifying we obtain the following 12 relations among the elements of $I_3(\Sigma)$:

\begin{eqnarray}
\ell_1\ell_3\ell_4-\ell_2\ell_3\ell_4&=&0\\
\ell_1\ell_3\ell_5-\ell_2\ell_3\ell_5&=&0\\
\ell_1\ell_4\ell_5-\ell_2\ell_4\ell_5&=&0\\
\ell_1\ell_2\ell_4-\ell_2\ell_3\ell_4&=&0\\
\ell_1\ell_2\ell_5-\ell_2\ell_3\ell_5&=&0\\
\ell_1\ell_4\ell_5-\ell_3\ell_4\ell_5&=&0\\
\ell_1\ell_2\ell_4-\ell_1\ell_3\ell_4&=&0\\
\ell_1\ell_2\ell_5-\ell_1\ell_3\ell_5&=&0\\
\ell_2\ell_4\ell_5-\ell_3\ell_4\ell_5&=&0\\
\ell_1\ell_2\ell_4-\ell_1\ell_2\ell_5&=&0\\
\ell_1\ell_3\ell_4-\ell_1\ell_3\ell_5&=&0\\
\ell_2\ell_3\ell_4-\ell_2\ell_3\ell_5&=&0.
\end{eqnarray}

Because $(\ell_2-\ell_3)=(\ell_1-\ell_3)-(\ell_1-\ell_2)$, relations $(7)$, $(8)$, and $(9)$ can be obtained immediately from the first six. With $(11)=(10)-(7)+(8)$ and $(12)=(10)-(4)+(5)$, we obtain that $(1),\ldots,(6),(10)$ are linearly independent, and so $N_{3,\Sigma}=7$. Immediately, by Theorem \ref{thm_kernel}~(iii), we obtain that $$b_1(3,\Sigma)={{5}\choose{2}}-N_{3,\Sigma}=10-7=3.$$

If we order lexicographically the columns of the $7\times {{5}\choose{2}}$ matrix $\mathcal M$ (i.e., $(1,2), (1,3),\ldots,(4,5)$) we obtain

$$\mathcal M=\left[\begin{array}{rrrrrrrrrr}0&0&0&-1&0&0&1&0&0&0\\
0&0&-1&0&0&1&0&0&0&0\\
0&-1&0&0&1&0&0&0&0&0\\
0&0&0&-1&0&0&0&0&1&0\\
0&0&-1&0&0&0&0&1&0&0\\
-1&0&0&0&1&0&0&0&0&0\\
0&0&0&0&0&0&0&-1&1&0\end{array}\right].$$ Observe that the $(4,5)$ column corresponding to the circuit $\{4,5\}$ is a column of all zeros, which agrees with the discussion in Remark \ref{rem_kernelPsi}~(a).

Let $i\geq 1$ and suppose $(g_1,\ldots,g_7)\in (\ker\Psi_{3,\Sigma})_i$. Then we obtain:
\begin{eqnarray}
-g_6\in\langle\ell_1,\ell_2\rangle=\langle x\rangle\nonumber\\
-g_3\in\langle\ell_1,\ell_3\rangle=\langle x\rangle\nonumber\\
-g_2-g_5\in\langle\ell_1,\ell_4\rangle=\langle x,y\rangle\nonumber\\
-g_1-g_4\in\langle\ell_1,\ell_5\rangle=\langle x,y\rangle\nonumber\\
g_3+g_6\in\langle\ell_2,\ell_3\rangle=\langle x\rangle\nonumber\\
g_2\in\langle\ell_2,\ell_4\rangle=\langle x,y\rangle\nonumber\\
g_1\in\langle\ell_2,\ell_5\rangle=\langle x,y\rangle\nonumber\\
g_5-g_7\in\langle\ell_3,\ell_4\rangle=\langle x,y\rangle\nonumber\\
g_4+g_7\in\langle\ell_3,\ell_5\rangle=\langle x,y\rangle.\nonumber
\end{eqnarray}

By Remark \ref{rem_kernelPsi}~(b), there is no restriction on $g_1,g_2,g_4,g_5,g_7\in R_i$, and $g_3,g_6\in\langle x\rangle_i$. So, for $i\geq 1$, $$\dim_{\mathbb K}(\ker\Psi_{3,\Sigma})_i=5(i+1)+2i=7i+5.$$

By Theorem \ref{thm_summary}~(5), $b_1(4,\Sigma)={\rm HF}(I_4(\Sigma),4)=2$, and so applying Theorem \ref{thm_kernel}~(iv) with $i=1$, we obtain
${\rm HF}(I_3(\Sigma),4)-2=4-7\cdot 2+12=2$, and so ${\rm HF}(I_3(\Sigma),4)=4$. But, from Lemma \ref{lem_HF} we have
$$b_2(3,\Sigma)=2b_1(3,\Sigma)-{\rm HF}(I_3(\Sigma),3+1)=2.$$

Similarly, since ${\rm HF}(I_4(\Sigma),5)=3$, applying Theorem \ref{thm_kernel}~(iv) with $i=2$ we obtain ${\rm HF}(I_3(\Sigma),5)=3+4-7\cdot 3+19=5$, and therefore $$b_3(3,\Sigma)=2b_2(3,\Sigma)-3b_1(3,\Sigma)+{\rm HF}(I_3(\Sigma),3+2)=0,$$ as expected.
\end{exm}

\begin{rem}\label{rem_revisit} Since the complex $(*)$ was inspired by the work in \cite{To3}, we find it appropriate to end with a remark where we recover the main result in this article, result that we mentioned also in part (iii) of Section \ref{hyperplanes}. Let $\Sigma=(\ell_1,\ldots,\ell_n)\subset R:=\mathbb K[x_1,\ldots,x_k]$ be a collection of linear forms defining a hyperplane arrangement of rank $k$. Let $a=n-2$. The goal is to use Theorem \ref{thm_kernel} to produce the Betti numbers for $I_{n-2}(\Sigma)$.

\medskip

$\bullet$ First, we determine $N_{n-2,\Sigma}$. Since $\Sigma$ defines a hyperplane arrangement, then all circuits will have length $\geq 3$. Let $1\leq s\leq n-(n-2)+1=3$ and look for circuits of length $s$ to build $\mathcal W(n-2,\Sigma)$. So $s=3$.

If $X_1,\ldots,X_{\gamma}$ are all distinct the flats of rank 2 of $\Sigma$, then for each $j=1,\ldots,\gamma$, $X_j$ will produce $\displaystyle {{|\Sigma_{X_j}|}\choose{3}}$ circuits of length $s=3$; recall the notations from Section \ref{hyperplanes} (iii). Also, each such circuit will produce $\displaystyle {{n-3}\choose{n-3-(n-2)+1}}={{n-3}\choose{0}}=1$ vectors in $\mathcal W(n-2,\Sigma)$.

Let $X\in L_2(\Sigma)$ and suppose $\Sigma_X=\{\ell_1,\ldots,\ell_{\alpha}\}, \alpha\geq 3$. For $1\leq u<v<w\leq \alpha$, consider the circuit $\{u,v,w\}$ with the corresponding dependency $D_{u,v,w}:\, c_{v,w}\ell_u+c_{u,w}\ell_v+c_{u,v}\ell_w=0$. This produces the only relation $$c_{v,w}F_{[n]\setminus\{v,w\}}+c_{u,w}F_{[n]\setminus\{u,w\}}+c_{u,v}F_{[n]\setminus\{u,v\}}=0,$$ and therefore the vector in $\mathcal W(n-2,\Sigma)$ $${\bf c}_{u,v,w}:=(0,\ldots,0,c_{u,v},0,\ldots,0,c_{u,w},0,\ldots,0,c_{v,w},0,\ldots,0),$$ where $c_{i,j}$ sits on the $(i,j)$-th entry, ordered by lexicographic order.

After a change of coordinates, we can assume that $\ell_i=x_1+\lambda_ix_2, i=1,\ldots,\alpha$, with $\lambda_1=0$ and $\lambda_i\neq\lambda_j, i\neq j$. Then $c_{v,w}=\lambda_v-\lambda_w$, $c_{u,w}=\lambda_w-\lambda_u$, and $c_{u,v}=\lambda_u-\lambda_v$, and if $u\geq 2$, we have the dependencies

$${\bf c}_{u,v,w}=\frac{1}{\lambda_u}[(\lambda_u-\lambda_w){\bf c}_{1,u,v}-(\lambda_u-\lambda_v){\bf c}_{1,v,w}].$$ This implies that $X$ contributes with $\displaystyle {{|\Sigma_{X}|-1}\choose{2}}$ linearly independent vectors in $\mathcal W(n-2,\Sigma)$.

Lastly, if ${\bf c}_{u,v,w}$ and ${\bf c}_{u',v',w'}$ corresponding to two distinct rank 2 flats $X$ and $X'$ have nonzero entries on the same place $(i,j)$, then $\{u,v,w\}\cap\{u',v',w'\}=\{i,j\}$ and so $X=X'=V(\ell_i,\ell_j)$; contradiction. So, with the convention that ${{2}\choose{3}}=0$, we have $$N_{n-2,\Sigma}=\sum_{X\in L_2(\Sigma)}{{|\Sigma_{X}|-1}\choose{2}}.$$

By Theorem \ref{thm_kernel}~(iii), $$b_1(n-2,\Sigma)=\underbrace{{{n}\choose{2}}}_{\alpha}-\underbrace{\sum_{X\in L_2(\Sigma)}{{|\Sigma_{X}|-1}\choose{2}}}_{\beta},$$ exactly what it is stated in Section \ref{hyperplanes}~(iii).

\medskip

$\bullet$ Next, we determine $\dim_{\mathbb K}(\ker\Psi_{n-2,\Sigma})_{i}$, for $i\geq 0$. Since each distinct $X$ of rank 2 produces its own direct summand in $\mathcal W(n-2,\Sigma)$, it is enough to see how $(\ker\Psi_{\alpha-2,\Sigma_X})_{i}, i\geq 1$, looks ``locally'' at $X$, where $\gamma=|\Sigma_X|\geq 3$. As before, we can assume $\ell_i=x_1+\lambda_ix_2, i=1,\ldots,\gamma$, with $\lambda_1=0$ and $\lambda_i\neq\lambda_j, i\neq j$.

For any $1\leq u<v\leq\gamma$, let $g_{u,v}\in R_i$ be such that $$\sum_{2\leq u<v\leq\gamma}\hat{g}_{u,v}{\bf c}_{1,u,v}=0\in \bigoplus_{2\leq u<v\leq\gamma}\frac{R}{\langle\ell_u,\ell_v\rangle}.$$

Since each ${\bf c}_{1,u,v}$ is built from the linear dependency $c_{1,u}\ell_v+c_{1,v}\ell_u+c_{u,v}\ell_1=0$ we obtain that for each $2\leq u\leq \gamma$, $$c_{1,u}(g_{2,u}+\cdots+g_{u-1,u}+g_{u,u+1}+\cdots+g_{u,\alpha})\in\langle\ell_1,\ell_u\rangle$$ and for any $1\leq u<v\leq \gamma$, $$c_{u,v}g_{u,v}\in\langle\ell_u,\ell_v\rangle.$$ But the first set of conditions are redundant: for example, since $g_{2,u}=A\ell_2+B\ell_u$, then $$c_{1,u}g_{2,u}=A(-c_{1,2}\ell_u-c_{2,u}\ell_1)+c_{1,u}B\ell_u\in\langle\ell_1,\ell_u\rangle.$$ So the dimension of this part of $(\ker\Psi_{n-2,\Sigma})_{i}$ is $$\sum_{2\leq u<v\leq \gamma}\dim_{\mathbb K}\langle \ell_u,\ell_v\rangle_i= {{\gamma-1}\choose{2}}\left(2{{k+i-2}\choose{i-1}}-{{k+i-3}\choose{i-2}}\right).$$ If we sum up we obtain that for $i\geq 0$,
$$\dim_{\mathbb K}(\ker\Psi_{n-2,\Sigma})_{i}=\beta\left(2{{k+i-2}\choose{i-1}}-{{k+i-3}\choose{i-2}}\right).$$

\medskip

$\bullet$ Now we compute $b_i(n-2,\Sigma)$, for $i\geq 2$.

Since $\displaystyle {\rm HF}(I_{n-1}(\Sigma),n-2+i-1)=n{{k+i-3}\choose{i-2}}-(n-1){{k+i-4}\choose{i-3}}$, Theorem \ref{thm_kernel}~(iv) gives
\begin{eqnarray}
{\rm HF}(I_{n-2}(\Sigma),n-2+i-1)&=&n{{k+i-3}\choose{i-2}}-(n-1){{k+i-4}\choose{i-3}}+\alpha{{k+i-4}\choose{i-1}}\nonumber\\
                           &-&\beta\left({{k+i-2}\choose{i-1}}-2{{k+i-3}\choose{i-2}}+{{k+i-4}\choose{i-3}}\right).\nonumber
\end{eqnarray}
\begin{itemize}
  \item[(a)] If $i=2$, then ${\rm HF}(I_{n-2}(\Sigma),n-2+1)=n+\alpha(k-2)-\beta(k-2)$, and so, from Lemma \ref{lem_HF} we have $$b_2(n-2,\Sigma)= k(\alpha-\beta)-n-\alpha(k-2)+\beta k-2\beta=2\alpha-n-2\beta,$$ exactly what we obtained in Section \ref{hyperplanes}~(iii).
  \item[(b)] If $i=3$, from Lemma \ref{lem_HF} we have $$b_3(n-2,\Sigma)=kb_2(n-2,\Sigma)-{{k+1}\choose{2}}b_1(n-2,\Sigma)+{\rm HF}(I_{n-2}(\Sigma),n-2+2),$$ and if we plug in we obtain the desired $$b_3(n-2,\Sigma)=\alpha-n-\beta+1.$$ One could obtain the same result by using part (c) below, and the fact that $b_1(n-2,\Sigma)-b_2(n-2,\Sigma)+b_3(n-2,\Sigma)=1$.
  \item[(c)] If $i=4$, because we know that ${\rm pdim}_R(R/I_{n-2}(\Sigma))=\min\{k,n-(n-2)+1\}=3$, as $k\geq 3$, we obtain that $b_4(n-2,\Sigma)=0$. We leave it to the reader to check that we get the same result if we use the formulas above.
\end{itemize}
\end{rem}


\bigskip

\renewcommand{\baselinestretch}{1.0}
\small\normalsize 

\bibliographystyle{amsalpha}

\end{document}